\documentclass[10pt,a4paper,twoside]{article}
\usepackage{amsmath, exscale, epsfig,amssymb,makeidx, mathrsfs}

\small
\normalsize
\usepackage[latin1]{inputenc}
\usepackage[T1]{fontenc}

\usepackage{amsfonts}

\usepackage{amsbsy}

\usepackage{amscd}
\usepackage{theorem}

\usepackage{epsf}

\usepackage{psfrag}
\usepackage{graphicx}

\usepackage{multicol}
\allowdisplaybreaks

\def\build#1_#2^#3{\mathrel{\mathop{\kern 0pt#1}\limits_{#2}^{#3}}}
\def\noi{{\noindent}}

\def\qed{$\hfill \square$}
\def\cq{$\hfill \square$}
\def\un{{\bf 1}}
\newcommand{\bbD}{\mathbb{D}}
\newcommand{\bD}{\mathbb{D}}
\newcommand{\bbE}{\mathbb{E}}
\newcommand{\bE}{{\bf E}}

\newcommand{\bm}{{\bf m}}

\newcommand{\bN}{\mathbb{N}}

\newcommand{\bbP}{\mathbb{P}}

\newcommand{\bP}{{\bf P}}
\newcommand{\bR}{\mathbb{R}}
\newcommand{\R}{\mathbb{R}}

\newcommand{\bbT}{\mathbb{T}}

\newcommand{\cA}{\mathcal{A}}
\newcommand{\cB}{\mathcal{B}}
\newcommand{\cC}{\mathcal{C}}
\newcommand{\cI}{\mathcal{I}}

\newcommand{\cF}{\mathcal{F}}
\newcommand{\cG}{\mathcal{G}}

\newcommand{\cL}{\mathcal{L}}
\newcommand{\cM}{\mathcal{M}}
\newcommand{\cN}{\mathcal{N}}

\newcommand{\cP}{\mathscr{P}}
\newcommand{\ccK}{\mathscr{K}}
\newcommand{\ccZ}{\mathscr{Z}}

\newcommand{\ccS}{\mathscr{S}}

\newcommand{\ccC}{\mathscr{C}}

\newcommand{\cT}{\mathcal{T}}

\newcommand{\sigA}{\sigma \! \mathcal{A}}
\def\dHp{{\rm d}^{_{ (p) }}_{^{{\bf H}}}}

\def\cJ{{\cal J}}
\def\cD{{\cal D}}

\def\varep{\varepsilon}

\def\be{\begin{equation}}
\def\ee{\end{equation}}
\def\ba{\begin{eqnarray*}}
\def\ea{\end{eqnarray*}}

\def\noi{\noindent}

\newcommand{\lgeo}{[{\!  } [}
\newcommand{\rgeo}{] {\!  } ]}

\def\cqfd{ \hfill $\blacksquare$ }
\def\btau{{\boldsymbol{\tau}}}
\def\inject{{\rm j}}

\newtheorem{theorem}{Theorem}[section]
\newtheorem{lemma}[theorem]{Lemma}

{\theorembodyfont{\rmfamily}\newtheorem{definition}{Definition}[section]}
{\theorembodyfont{\rmfamily}\newtheorem{notation}{Notation}[section]}
{\theorembodyfont{\rmfamily}\newtheorem{comment}{Comment}[section]}
{\theorembodyfont{\rmfamily}\newtheorem{remark}{Remark}[section]}
\begin{document}

\title{ {\bf THE EXACT PACKING MEASURE OF LEVY TREES}.}
\author{Thomas {\sc Duquesne}  
\thanks{Laboratoire de Probabilit\'es et Mod\`eles Al\'eatoires; 
Universit\'e Paris 6, Bo\^ite courrier 188, 4 place Jussieu, 75252 Paris Cedex 05, FRANCE. Email: thomas.duquesne@upmc.fr. This works benefited from the partial support of ANR A3, Projet BLAN-****} }

\vspace{2mm}
\date{\today} 

\maketitle

\begin{abstract} 
We study fine properties of L\'evy trees that are random compact metric spaces introduced by Le Gall and Le Jan in 1998 as the genealogy of continuous state branching processes. L\'evy trees are the scaling limits of Galton-Watson trees and they generalize Aldous's continuum random tree which corresponds to the Brownian case. In this paper we prove that L\'evy trees have always an exact packing measure: We explicitely compute the packing gauge function and we prove that the corresponding packing measure coincides with the mass measure up to a multiplicative constant. 

\medskip

\end{abstract}

\section{Introduction}
\label{introsec}

  L\'evy trees are random compact metric spaces introduced by Le Gall and Le Jan in \cite{LGLJ1}   as the genealogy of Continuous State Branching Processes (CSBP for short). The class of L\'evy trees comprehends Aldous's continuum random tree which corresponds to the Brownian case (see \cite{Al1, Al2}). L\'evy trees are the scaling limits of Galton-Watson trees (see \cite{DuLG} Chapter 2 and \cite{Du2}). Various geometric and distributional properties of L\'evy trees have been studied in \cite{DuLG2}, in Weill \cite{Weill} and in Abraham and Delmas \cite{AbDel09}. An alternative construction of L\'evy trees is discussed in \cite{DuWi1}. L\'evy trees have also been studied in connection with fragmentation processes: see Miermont \cite{Mier03, Mier05}, Haas and Miermont \cite{HaaMi},  Goldschmidt and Haas \cite{GolHaa} for the stable cases; see Abraham and Delmas \cite{AbDel08} for more general models.

   Fractal properties of L\'evy have been discussed in \cite{DuLG2} and \cite{DuLG3}: Hausdorff and packing dimensions of L\'evy trees are computed in \cite{DuLG2} and the exact Hausdorff measure of the continuum random tree is given in \cite{DuLG3}. As shown in \cite{Du10} (see also 
\cite{DuLG3}), there is no exact Hausdorff measure in the non-Brownian stable cases. The goal of this paper is to prove that L\'evy trees behave better with respect to packing measure. More precisely, we prove that  L\'evy trees have always an exact packing measure: We explicitely compute the packing gauge function and we prove that the corresponding packing measure coincides with the mass measure up to a multiplicative constant. 
   
\medskip   
  
Before stating the main results of the paper, let us recall basic facts on continuous state branching processes (CSBP) and on L\'evy trees. CSBPs are time- and space-continuous analogues of Galton-Watson Markov chains. They have been introduced by Jirina \cite{Ji} and Lamperti \cite{La2} as the 
$[0, \infty]$-valued Feller processes that are absorbed in states $0$ and $ \infty $ and whose kernel semi-group $(p_t (x,dy); x \in [0, \infty], t \in [0, \infty))$ enjoys the {\it branching property}: $p_t (x, \cdot )*p_t (x^\prime,\cdot ) = p_t(x+ x^\prime , \cdot)$, for every $x,x^\prime\in [0, \infty]$ and every $t \in [0, \infty)$ (here $*$ stands for the convolution product).  
Let  $Z= (Z_t, t \geq 0)$ be a CSBP with initial state $x\in (0, \infty)$ that is defined on a probability space $(\Omega, \cF, \bP)$. {\it We shall restrict our attention to CSBPs that get almost surely extinct in finite time}. Namely, 
\begin{equation}
\label{extinctexplici}
\bP (\exists t \geq 0: \; Z_t= 0)= 1\; .
\end{equation} 
Then, Silverstein \cite{Sil68} proves that the kernel semigroup of $Z$ is characterised by a function $\psi: [0, \infty) \rightarrow [0, \infty)$ as follows: For any $\lambda, s,t \geq 0$, one has $\bE [ \exp (-\lambda Z_{t+s}) | Z_s]= \exp (-Z_s u(t, \lambda))$, where $u(t,\lambda)$ is the unique nonnegative solution of $\partial u (t , \lambda) /
 \partial t =-\psi (u(t, \lambda))$ with $u(0, \lambda)= \lambda$. This equation can be rewritten in  the following integral form. 
 \begin{equation}
\label{integral}
\int_{u(t, \lambda)}^\lambda \frac{du}{\psi (u)} = t \; , \quad t, \lambda \geq 0 . 
\end{equation}
The function $\psi$ is called {\it the branching mechanism} of the CSBP. Under Assumption (\ref{extinctexplici}), $\psi$ is necessarily of the following L\'evy-Khintchine form
\begin{equation}
\label{LevyKhin}
 \psi (\lambda )= \alpha \lambda + \beta \lambda^2 + 
\int_{(0, \infty)} (e^{-\lambda r} -1+\lambda r)\,\pi (dr)  \; ,  \quad \lambda \geq 0,
\end{equation}
where $\alpha \in [ 0, \infty)$ is the {\it drift coefficient}, $\beta \in [ 0, \infty)$ is the {\it Brownian coefficient} and $\pi$ is the {\it L\'evy measure} that satisfies $\int_{{(0, \infty)}} (r\wedge r^2)\,\pi (dr) < \infty$.  Moreover, $\psi$ has to satisfy the following condition: 
\begin{equation}
\label{extinction}
\int^\infty \frac{du}{\psi (u)} < \infty . 
\end{equation}
More precisely, the set of branching mechanisms of CSBPs that satisfy (\ref{extinctexplici})  is exactly the set of functions $\psi$ of the form (\ref{LevyKhin}) that satisfy (\ref{extinction}) (see  Bingham \cite{Bi2} for more details on CSBPs).

   Let us introduce the formalism developed in \cite{DuLG2} where L\'evy trees are viewed as 
random variables taking their values in the space of compact rooted $\bR$-trees.
Informally, a {\it $\bR$-tree} is a metric space $(\cT,d)$ such that for any two points
$\sigma$ and $\sigma^\prime $
in $\cT$, there is a unique arc with endpoints $\sigma$ and $\sigma^\prime$ and this arc is isometric to a compact interval of the real line. A {\it rooted $\bR$-tree} is a $\bR$-tree with a distinguished vertex $\rho$ called the {\it root}. We say that two rooted $\bR$-trees are {\it equivalent} if there is a
root-preserving isometry that maps one onto the other. Instead of considering all compact rooted $\bR$-trees, we introduce the set $\bbT$ of equivalence classes of compact rooted
$\R$-trees. Evans, Pitman and Winter \cite{EvPitWin} prove that $\bbT$ equipped
with the Gromov-Hausdorff distance \cite{Gro} is a Polish space. As proved in \cite{DuLG2}, with any branching mechanism $\psi$ of the form (\ref{LevyKhin}) that satisfies (\ref{extinction}), 
one can associate a sigma-finite measure $\Theta_\psi$ on $\bbT$ that is called the {\it "distribution'' of the L\'evy tree with branching mechanism $\psi$}. Although $\Theta_\psi$ is an infinite measure, the following holds true: Set $\Gamma (\cT) = \sup_{\sigma \in \cT} d(\rho, \cT)$, that is the total height of $\cT$; Then, for any $a \in (0, \infty)$, one has 
\begin{equation}
\label{intrinvv}
 v(a):=\Theta_\psi ( \Gamma (\cT)> a ) < \infty \; , 
\end{equation}
where the function $v: (0, \infty) \rightarrow (0, \infty)$ is determined by  $\int_{{v(a)}}^{\infty} \psi(u)^{-1}du=a$.

  L\'evy trees enjoy the so-called {\it branching property}, that obviously holds true for Galton-Watson trees: For every $a>0$, under the probability measure $\Theta_\psi (\, \cdot \, | \, \Gamma (\cT)>a)$ and conditionally given the part of $\cT$ below level $a$, the subtrees above level $a$ are distributed as the atoms of a Poisson point measure whose intensity is a random multiple of $\Theta_\psi $. The random factor is the total mass of the $a$-local time measure that is defined below. Let us mention that Weill \cite{Weill} shows that the branching property characterizes L\'evy trees.

$\Theta_\psi$ is approximated by Galton-Watson trees as follows. Let $(\xi_p)_{p \geq 0}$ be a sequence of probability distributions on the set of nonnegative integers $\bN$. We first assume that  $\sum_{{k\in \bN}} k\xi_p (k) \leq 1$, for any $p$ and that the $\xi_p$s are in the domain of attraction of an infinitely divisible distribution with Laplace exponent $\psi$. More precisely, let $Y$ be a real valued random variable such that $\log \bE [\exp(-\lambda Y) ]= \psi (\lambda)$, for any $\lambda \in [0, \infty)$. For any $p$, let $(J^{_{(p)}}_k)_{ k \geq 0}$ be an i.i.d.$\;$sequence of r.v.$\;$with law $\xi_p$. We assume there exists an increasing sequence $(b_p)_{p \geq 0}$
of positive integers such that $\frac{p}{b_p}(J^{_{(p)}}_1 +\cdots+J^{_{(p)}}_{b_p} -b_p)$
converges in distribution to $Y$. For every $p$, denote by $\btau_p$
a Galton-Watson tree with offspring distribution $\xi_p$ that can be viewed as a random rooted 
$\R$-tree $(\btau_p, d_p , \rho_p)$ by affecting unit length to each edge. Thus, $(\btau_p, \frac{_1}{^p}d_p , \rho_p)$ is the tree $\btau_p$ whose edges are rescaled by a factor $1/p$ and we simply denote it by $ \frac{_1}{^p} \btau_p$. We furthermore assume that for any $a \in (0, \infty)$, one has $\liminf_p \bP ( \Gamma (\btau_p)  \leq p \,  a  )^{b_p/p} >0$. Roughly speaking, this assumption ensures that 
$\frac{_1}{^p}\Gamma (\btau_p)$ has a non trivial limit in law. Under these assumptions, Theorem 4.1 \cite{DuLG2} asserts that for any $a \in (0 , \infty)$, {\it the law of $ \frac{_1}{^p} \btau_p$ under $\bP (\, \cdot \, | \, \frac{_1}{^p}\Gamma (\btau_p)  > a)$ converge weakly on $\bbT$ 
to $\Theta_\psi (\, \cdot \, | \, \Gamma (\cT)  > a)$, when $p$ goes to $\infty$}.

There are two important kinds of measures on $\psi$-L\'evy trees. For every $a>0$, let us set $\cT(a):=\{\sigma\in\cT :d(\rho,\sigma)=a\}$ that is the {\it $a$-level set of $\cT$}. 
Then, we define a measure $\ell^a$ on $\cT(a)$ as follows: For every $\varepsilon>0$, write
$\cT_\varepsilon(a)$ for the finite subset of $\cT(a)$ consisting of
those vertices which have descendants at level $a+\varepsilon$. Then, $\Theta_\psi$-a.e.$\;$for every
bounded continuous function $f$ on $\cT$, we have 
\begin{equation}
\label{geolocapprox}
\langle \ell^a,  f \rangle=\lim_{\varepsilon\downarrow 0} \; \frac{_1}{^{v(\varepsilon)}} \!\!\! \sum_{\sigma\in \cT_\varepsilon(a)} \!\!\!  f(\sigma), 
\end{equation}
where $v$ is defined by (\ref{intrinvv}). The finite measure $\ell^a$ on $\cT(a)$ is called the {\it $a$-local time measure of $\cT$}. We refer to Section 4.2 \cite{DuLG2} for the construction and the main properties of local time measures. Theorem 4.3 \cite{DuLG2} ensures that one can choose a modification of $(\ell^a, a \geq 0)$ such that $a \mapsto \ell^a$ is $\Theta_\psi$-a.e.$\;$cadlag with respect to the weak topology on finite measures on $\cT$. We also define the {\it mass measure} $\bm$ on the tree $\cT$ by 
\begin{equation}
\label{unifmeas}
\bm=\int_0^\infty da\,\ell^a\,.
\end{equation}
The topological support of $\bm$ is $\cT$ and $\bm$ is in some sense the most spread out measure on $\cT$.  Note that the definitions of the local time measures and of the mass measure $\bm$ only involve the metric properties of $\cT$. 

Let us recall from \cite{DuLG2} results concerning the Hausdorff and the packing dimensions of $\psi$-L\'evy trees:  Let $\gamma$ (resp. $\eta$) be the {\it lower} (resp. {\it upper}) {\it exponent of $\psi$ at infinity}. Namely, $\gamma = \sup \{ c \geq 0 : \; \lim_{ \infty} \psi (\lambda) \lambda^{-c} = \infty  \} $ and 
$\eta= \inf \{ c \geq 0  : \; \lim_{ \infty} \psi (\lambda) \lambda^{-c} = 0 \} $. Since $\psi$ is of the form (\ref{LevyKhin}), one clearly has $1\leq   \gamma \leq \eta \leq 2 $.  
Theorem 5.5 \cite{DuLG2} asserts that if $\gamma >1$, then $\Theta_\psi$-a.e.$\, \cT$ has Hausdorff dimension $\eta/ (\eta-1)$ and packing dimension $\gamma/ (\gamma-1)$. In this paper we discuss finer results concerning the exact packing measure of L\'evy trees.  
Packing measures have been introduced by Taylor and Tricot in \cite{TaTr}. Though their construction is done in Euclidian spaces, it extends to metric spaces, and more specifically to L\'evy trees. More precisely, for any $\sigma \in \cT$ and any $r \in [0, \infty)$, we denote by $\bar{B} (\sigma , r)$ (resp. $B(\sigma, r)$) the closed (resp. open) ball of $\cT$ with center $\sigma$ and radius $r$. Let $A \subset \cT$ and $\varepsilon \in (0,\infty  )$. A {\it $\varepsilon$-packing of $A$} is a countable collection of pairwise disjoint closed balls $\bar{B}(x_n , r_n) $,  $n \geq 0$,  such that $x_n \in A$  and $r_n \leq \varepsilon$. 
We restrict our attention to packing measures associated with {\it regular gauge functions} in the following sense: a function $g: (0, r_0) \rightarrow (0, \infty)$ is a regular gauge function if it is continuous, non decreasing, if $\lim_{0+} g= 0$  and if the following holds true 
\begin{equation}
\label{doubling}
\exists \, C >1 \; : \quad g(2r) \leq C g(r) \;, \quad r \in (0, r_0/2) . 
\end{equation}
Such property shall be refered as to a  {\it $C$-doubling condition}. The $g$-packing measure on $(\cT, d)$ is then defined as follows. 
For any $\varepsilon \in (0, \infty)$, we first set: 
\begin{equation}
\label{prepremeadef}
\cP^{(\varepsilon)}_g (A)=  \sup \Big\{ \sum_{^{n \geq 0}}  g(r_n) ;  
\;  \textrm{ $(\bar{B}(x_n , r_n) )_{ n \geq 0}$ is a  $\varepsilon$-packing of $A$}  \Big\} \; .
\end{equation}
 The {\it $g$-packing pre-measure of $A$} is then defined by 
\begin{equation}
\label{premeadef}
 \cP^*_g (A)=\lim_{\varepsilon \downarrow 0} \cP^{(\varepsilon)}_g (A) 
\end{equation}
and we define the $g$-{\it packing outer measure of $A$} as 
\begin{equation}
\label{packdef}
\cP_g (A) = \inf \Big\{ \sum_{^{n \geq 0}} \cP^*_g (E_n) ; \;  A \subset \bigcup_{^{n \geq 0}} E_n  \Big\} \; .
\end{equation}
As in Euclidian spaces, $\cP_g$ is a Borel regular metric outer measure (see Section \ref{packingsec} for more details). The original definition of packing measures \cite{TaTr} makes use, as set function, of the diameter of open ball packing instead of the radius of closed ball packing. As pointed out by H. Haase \cite{Haase}, diameter-type packing measures may be not Borel regular: H. Joyce \cite{Joyce99} provides 
an explicit example where this problem occurs. In our setting, we don't face such problem and our 
results hold true for diameter-type and radius-type packing measures as well, thanks to specific 
properties of compact real trees (see \cite{Du10} for more details).

  We state below that  the $\psi$-L\'evy tree $\cT$ has an exact packing measure with respect to a gauge function that is defined as follows: Denote by $\psi^\prime$ the derivative of $\psi$ and note that $\psi^\prime (0)= \alpha$. We also denote by $\psi^{-1}$ the reciprocal of $\psi$. Since $\psi$ is of the form (\ref{LevyKhin}) and since it satisfies (\ref{extinction}), $\psi^\prime$ and $\psi^{-1}$ both tend 
to $\infty$ (see Section \ref{notadefsec} for more details). We then set $\varphi = \psi^\prime \circ \psi^{-1} $ and we denote by $\varphi^{-1} $ its reciprocal function that is defined from $[\alpha, \infty)$ to $[0, \infty)$. We then define the $\psi$-gauge function $g$ by 
\begin{equation} 
\label{gaugedef}    
 g (r) := \frac{\log \log\frac{1}{r}}{\varphi^{-1} \!\! \left( \frac{1}{r}\log \log \frac{1}{r} \right)} \; , \; r\in (0, r_0) 
\end{equation}
where $r_0$ stands for $\min (\alpha^{-1} , e^{-e})$ (with the convention $\alpha^{-1}= \infty$ if $\alpha= 0$). 

\begin{theorem}
\label{packingLevytreeth} Let $\psi$ be a branching mechanism of the form (\ref{LevyKhin}). We assume that the function $g$ that is derived from $\psi$ by (\ref{gaugedef}) satisfies a doubling condition (\ref{doubling}). Let $(\cT, d, \rho)$ be the $\psi$-L\'evy tree under its excursion measure $\Theta_\psi$. Then, there exists a constant $c_\psi \in (0, \infty)$, that only depends on $\psi$ such that 
$$ \Theta_\psi \,^{_{_-}} {\rm a.e.} \qquad c_\psi \, \cP_g = \bm \; .$$
\end{theorem}
Although it is possible to define packing measures associated with gauge functions that don't satisfy a doubling condition (see Edgar \cite{Edgar07}), it obviously leads to technical complications. 
To assume that the gauge function $g$ (given by (\ref{gaugedef})) satisfies a doubling condition (\ref{doubling}), is in some sense the minimal hypothesis on $\psi$ that is required to stay in the standard framework of packing measure theory.  Let us briefly discuss this assumption: Lemma \ref{doublinglemma} Section \ref{exponentsec} asserts that {\it $g$ given by (\ref{gaugedef}) satisfies a doubling condition (\ref{doubling}) iff $\delta >1$, where $\delta$ stands for the following exponent:}
\begin{equation}
\label{doublingexponent}
 \delta := \sup \{ c \geq 0 : \, \exists Q \in (0, \infty)\; {\rm s.t.} \quad  Q . \psi (u) u^{-c}  \!\! \leq \psi (v) v^{-c}  \, , \, 1 \leq u \leq v \, \}.
 \end{equation}
We obviously have $1\leq \delta \leq  \gamma \leq \eta \leq 2 $, where $\gamma$ and $\eta$ are the lower and the upper exponents of $\psi$ at $\infty$. As already mentioned, if $\gamma >1$, then the packing dimension of $\cT$ is $\Theta_\psi$-a.e.$\;$equal to $\gamma/ (\gamma -1)$. So, $\gamma >1$ may look as a more natural assumption for Theorem \ref{packingLevytreeth} to be true. However Lemma \ref{contreex} Section \ref{exponentsec} shows that for any $c \in (1, 2]$, there exists a branching mechanism $\psi$ that is of the form (\ref{LevyKhin}), that satisfies (\ref{extinction}) and such that $\gamma = \eta = c$ but $\delta = 1$. Thus, $\delta >1$ is a more restrictive assumption than $\gamma >1$. Let us first mention that $\delta >1$ implies (\ref{extinction}), which is therefore not explicitly assumed in 
Theorem \ref{packingLevytreeth}. Let us also mention that if $\psi$ is regularly varying at $\infty$, then $\delta = \gamma= \eta$. 

\medskip 

  One important argument of the proof of Theorem \ref{packingLevytreeth} is the following result that gives the lower density of $\bm$ for typical points. 
\begin{theorem}
\label{densityth} Let $\psi$ be a branching mechanism of the form (\ref{LevyKhin}). We assume that the function $g$ that is derived from $\psi$ by (\ref{gaugedef}) satisfies a doubling condition (\ref{doubling}). Let $(\cT, d, \rho)$ be the $\psi$-L\'evy tree under its excursion measure $\Theta_\psi$. Then, there exists a constant $C_\psi \in (0, \infty)$, that only depends on $\psi$ such that $\Theta_\psi$-a.e.$\;$for $\bm$-almost all $\sigma$, one has 
$$ \liminf_{r\rightarrow 0} 
\frac{\bm \big( B(\sigma , r) \big)}{g(r)} = C_\psi . $$
\end{theorem}

\begin{comment}
\label{afterthcom} {\bf (a)} Let us consider the stable L\'evy trees. Namely, the branching mechanism $\psi$ is of the form 
$\psi (\lambda)= \lambda^\gamma $, with $\gamma \in (1, 2]$. The packing gauge function can be taken as $g(r)=r^{_{\gamma/ (\gamma-1)}} / ( \log \log 1/r )^{_{1/ (\gamma-1)}}$ and (\ref{doubling}) is obviously satisfied. Let us mention that $\gamma$-stables trees enjoy the following scaling property: For any $C \in (0, \infty)$, the law of $(\cT , C d, \rho)$ under $\Theta_\psi$ is $C^{_{1/(\gamma -1)}}\Theta_\psi$. The scaling property makes possible to condition $\cT$ to have total mass equal to $1$ (if $\gamma =2$,  then $ \Theta_\psi (\, \cdot \, | \, \bm (\cT)= 1 )$ is the law of the continuum random tree). Easy arguments imply that Theorem \ref{packingLevytreeth} and Theorem \ref{densityth} hold true under $ \Theta_\psi (\, \cdot \, | \, \bm (\cT)= 1 )$. 

In the Brownian case $\gamma = 2$, the packing gauge function can be taken as $g(r)=r^{2}/ \log \log 1/r $. As shown in \cite{DuLG3}, there exists a constant $c^\prime_\psi$ such that $\Theta_\psi$-a.e.$\;$ $c^\prime_\psi  {\mathscr H}_h = \bm$, where  
${\mathscr H}_h$ stands for the $h$-Hausdorff measure with $h(r)=  r^{2} \log \log 1/r $. Note that $g$ and $h$ are resp. the packing and the Hausdorff gauge functions for Brownian motion in dimensions $d \geq 3$. Theorem 1.10 \cite{Du10} asserts that non-Brownian $\gamma$-stable L\'evy trees (i.e. $\gamma \in (1, 2)$) have no exact Hausdorff measure with regularly varying gauge function. 
Thus, Theorem \ref{packingLevytreeth} shows that  L\'evy trees behave better with respect to packing measures than with respect to Hausdorff measures.

\smallskip

\noi 
{\bf (b)} Although Theorem \ref{packingLevytreeth} asserts that whole stable L\'evy trees have an exact packing measure, the level sets of stable L\'evy trees have no exact packing measure, even in the Brownian case, as shown by Theorem 1.1 \cite{Du10}.

\smallskip

\noi
{\bf (c)} Since $(\cT, d, \rho)$ is not Euclidian, the constant $C_\psi$ in Theorem \ref{densityth} may be distinct from the constant $c_\psi$ in Theorem \ref{packingLevytreeth}. This causes technicalities in the proof of Theorem \ref{packingLevytreeth}. 
Let us mention that $C_{\psi}= \gamma -1$, when $\psi (\lambda)= \lambda^\gamma$, $\gamma \in (1, 2]$.

\smallskip

\noi
{\bf (d)} Let $(S_r , r \geq 0)$ be a subordinator with Laplace exponent $\varphi= \psi^\prime \circ \psi^{-1}$. If $g$ (given by (\ref{gaugedef})) satisfies a doubling condition (\ref{doubling}), Fristedt and Pruitt in \cite{FriPru71} prove that there exists a constant $K_\varphi \in (0, \infty)$ that only depends on $\varphi$ and such that $\liminf_0 S_r/g(r)= K_\varphi$ a.s. Let us mention that we rely on this result for the proof of the upper bound in Theorem \ref{densityth}. \cq 
\end{comment}

The paper is organised as follows. In Section \ref{packingsec} we recall basic properties of packing  measures in metric spaces and a comparison result. In Section \ref{Levytreesec}, we introduce the height processes, the L\'evy trees and a key decomposition of L\'evy trees according to the ancestral line of a randomly chosen vertex. This decomposition plays an important role in the proof of Theorem \ref{densityth}. In Section \ref{exponentsec}, we prove Lemma \ref{doublinglemma} that shows that the gauge function $g$ satisfies a doubling condition iff the exponent $\delta$ given by (\ref{doublingexponent}) is strictly larger that $1$. In Section \ref{estimsec} we prove various estimates that are used in the proof sections. Section \ref{proofdensthsec} and \ref{proofpackingth} are devoted to the proof of resp.Theorem \ref{densityth} and Theorem \ref{packingLevytreeth}. 

\section{Notation, definitions and preliminary results.}
\label{notadefsec}
\subsection{Packing measures on metric spaces.}
\label{packingsec}
As already mentioned, we restrict our attention to continuous increasing gauge functions that satisfy a doubling condition as defined by (\ref{doubling}).  
Let $(\cT, d)$ be an uncountable complete and separable metric space. Let us fix a  regular gauge function $g$. Recall from (\ref{premeadef}) the definition of the $g$-packing pre-measure $\cP^*_g$. The $g$-packing pre-measure is non decreasing with respect to inclusion, it is sub-additive and it is a metric set function. Namely, if $A$ and $B$ are non-empty subsets of $\cT$ and if $\inf_{\sigma \in A , \sigma^\prime \in B} d(\sigma, \sigma^\prime) > 0 $, then  $\cP^*_g (A \cup B) =\cP_g^* (A) + \cP^*_g (B)$. Moreover $\cP_g^*$ has the following property. For any $A \subset \cT$, denote by $\bar{A}$ the closure of $A$. Then, we have 
\begin{equation}
\label{closprepack}
\cP_g^* (A)= \cP_g^* ( \bar{A} ) \; .
\end{equation}
Recall from (\ref{packdef}) the definition of the $g$-packing outer measure $\cP_g$.  As proved in \cite{Edgar07} Section 5, $\cP_g$  is a metric Borel regular outer measure satisfying the following properties.  
\begin{itemize}
\item{{\footnotesize{\bf Pack(1)}}} For any $A \subset \cT$, $\cP_g (A) \leq \cP_g^* (A)$. 
\item{{\footnotesize{\bf Pack(2)}}} If $A$ is $\cP_g$-measurable and such that $0<\cP_g (A)< \infty$, 
then for any $\varepsilon >0$, there exists a closed set $F \subset A$ such that $ \cP_g (A) \leq \cP_g (F ) + \varepsilon $. 
\item{{\footnotesize{\bf Pack(3)}}} $ \cP_g ( A) \!= \! \inf \big\{  \sup_{n \geq 0} \cP^*_g (A_n)  ;  \, A_n \! \subset \! A_{n+1}  \; {\rm and} \;  \bigcup_{^{n \geq 0}} A_n = A  \big\}$, for any $A\subset \cT$.
\end{itemize}
\noi
We shall also use the following comparison Lemma.  
\begin{lemma}
\label{genedens} (Taylor and Tricot \cite{TaTr}  Theorem 5.4, Edgar \cite{Edgar07} Theorem 5.9 ). 
Let $g $ be a regular gauge function that satisfies a $C$-doubling condition. Then, for any finite Borel measure $\mu$ on $\cT$ and for any Borel subset $A$ of $\cT$, the following holds true.
\begin{itemize} 
\item[(i)] If  $ \liminf_{r \rightarrow 0} \frac{\mu (B(\sigma, r))}{g(r )} \leq 1 $ for any $\sigma \in A$, then $  
\cP_g (A)  \geq C^{-2} \mu (A) $. 
\item[(ii)] If  $ \liminf_{r \rightarrow 0} \frac{\mu (B(\sigma, r))}{g(r)} \geq 1 $ for any $\sigma \in A$, then 
$  \cP_g (A)  \leq    \mu (A) $. 
\end{itemize}
\end{lemma}

\subsection{Height processes and L\'evy trees.}
\label{Levytreesec}
In this section we recall (mostly from \cite{DuLG} and \cite{DuLG2}) various results concerning height processes and L\'evy trees for further use in Section \ref{estimsec}, Section \ref{proofdensthsec} and Section \ref{proofpackingth}. 

\smallskip

\noi
{\bf The height process.} Recall that $\psi$ stands for a branching mechanism of the form (\ref{LevyKhin}). {\it We always assume that $\psi $ satisfies (\ref{extinction})}. It is convenient to work on the canonical space $\bD ([0, \infty), \bR )$ of cadlag paths equipped with Skorohod topology and the corresponding Borel sigma-field. We denote by  $X= (X_t , t \geq 0)$ the canonical process and by $\bbP$ the distribution of the spectrally positive L\'evy processes with Laplace exponent $\psi$. Namely, $\bbE [\exp (-\lambda X_t)  ] = \exp (t \psi (\lambda)\, )$, $\lambda, t \geq 0$. Note that the specific form of $\psi$ implies that $X_t$ is integrable and that $\bbE [X_t]= -\alpha t$. This easily entails that $X$ does not drift to $\infty$. Moreover (\ref{extinction}) implies that either $\beta >0$ or $\int_{^{(0, 1)}} r \pi (dr) = \infty$. It entails that $\bbP$-a.s.$\, X$ has unbounded variation sample paths (see Bertoin \cite{Be} Chapter VII Corollary 5 (iii)).

 As shown by Le Gall and Le Jan \cite{LGLJ1} (see also \cite{DuLG} Chapter 1), there exists a {\it continuous process} $H= (H_t , t \geq 0)$ such that for any $t \in [0, \infty)$, the following limit holds true in $\bbP$-probability 
\begin{equation}
\label{Hlimit}
H_t:=\lim_{\varepsilon\to 0} \frac{1}{\varep}\int_0^t {\bf 1}_{\{I^s_t<X_s<I^s_t+\varep\}}\,ds . 
\end{equation}
Here $I^s_t$ stands for $\inf_{s\leq r\leq t} X_r$. We shall use the notation $I_t= I^0_t= \inf_{0\leq r\leq t} 
X_r$, for the infimum of $X$. The process $H= (H_t, t \geq 0)$ is called the $\psi$-{\it height process}. As we see below, $H$ provides a way to explore the genealogy of CSBPs. We refer to Le Gall and Le Jan \cite{LGLJ1} for an explanation of (\ref{Hlimit}) in terms of discrete processes.

\smallskip

\noi
{\bf Excursions of the height process.} When $\psi$ is of the form $\psi (\lambda) = \beta \lambda^2$, $X$ is distributed as a Brownian motion and (\ref{Hlimit}) easily implies that $H$ is proportional to $X-I$, which is distributed as a reflected Brownian motion. In the general cases, $H$ is neither a Markov process nor a martingale. However it is possible to develop an excursion theory for $H$ as follows.  Recall that (\ref{extinction}) entails that $X$ has unbounded variation sample paths. Basic results on fluctuation theory (see Bertoin \cite{Be} Chapters VI.1 and VII.1) entail that $X-I$ is a strong Markov process in $[0, \infty)$ and that $0$ is regular for 
$(0, \infty)$ and recurrent with respect to this Markov process. Moreover, $-I$ is a local time at $0$ for $X-I$ (see Bertoin \cite{Be} Theorem VII.1). We denote by $N$ the corresponding excursion 
measure of $X-I$ above $0$. We denote by $(a_j, b_j)$, $j\in  \cI$, the excursion intervals of $X-I$ above $0$, and by $X^j = X_{(a_j + \cdot )\wedge b_j}-I_{a_j}$, $j\in \cI$, the corresponding excursions. Then, $\sum_{j\in \cI} \delta_{(-I_{a_j}, X^j)}$ is a Poisson point measure on $[0, \infty)\times \bbD([0, \infty), \bR)$ with intensity 
$dx \otimes N (dX)$. First observe that under $\bbP $, the value of $H_t$ only depends on the excursion of $X-I$ straddling $t$. Next note that  $\bigcup_{^{j\in \cI}} (a_j, b_j)= \{ t \geq 0: H_t >0 \}$. 
This allows to define 
the height process under $N$ as a certain measurable function $H(X)$ of $X$. 
See \cite{DuLG} Chapter 1, for more details.

\begin{notation}
\label{excuheightdef} Let $\cC^0$ be the space of the continuous functions from $[0, \infty)$ to $\bR$ equipped with the topology of the uniform convergence on every compact subsets of $[0, \infty)$ that makes it a Polish space. We shall denote by $\mathscr{C}$ the set of functions $h \in \cC^0$ with {\it compact support}. For any $h\in \ccC$, we set $\zeta(h)= \sup \{ t \in [0, \infty): h(t) \neq 0 \}$, with the convention $\sup \emptyset = 0$. If $h \in \cC^0 \backslash \ccC$, then $\zeta (h)= \infty$. 
By convenience, we denote by $H= (H_t , \geq 0)$ {\it the canonical process on $\cC^0$} and we call $\zeta = \zeta (H)$ the {\it lifetime of $H$}. 
{\it We slightly abuse notation by denoting by $N(dH)$ the "distribution" of the height process $H(X)$ associated with $X$ under the excursion measure $N (dX)$}. \cq 
\end{notation}

Note that $N$-a.e.$\; \zeta < \infty$, $H_0=H_{\zeta}=0$ and $H_t >0$ for any $t \in (0, \zeta)$. 
We now recall the Poisson decomposition of the height process $H(X)$ associated with $X$ under $\bbP$.  Recall that the intervals $(a_j, b_j)$, $j \in \cI$, are the open connected components of the set $\{ t \geq 0 : H_t >0\}$. For any $j\in  \cI$, we set $H^j =H_{(a_j+\cdot )\wedge b_j} $.Then, under $\bbP$, the point measure 
\begin{equation}
\label{Poissheight}
\sum_{j\in \cI} \delta_{(-I_{a_j},H^j)}
\end{equation}
is distributed as a Poisson point measure on $[0, \infty)\times\cC^0$ with intensity $dx \otimes N(dH)$. Note that under $N$, $X$ and $H$ have the same lifetime and recall that basic results of fluctuation theory entail
\begin{equation}
\label{lifetimeexc}
N \left( 1-e^{-\lambda \zeta} \, \right)= \psi^{-1} (\lambda ) \; , \quad \lambda \geq 0 . 
\end{equation}

\smallskip

\noi
{\bf Local times of the height process.} We recall  from \cite{DuLG} Chapter 1 Section 1.3 the following result: There exists a jointly measurable process $(L^a_s, a, s\geq 0)$ such that $\bbP$-a.s.$\,$ for any $a\geq 0 $, $s \mapsto L^a_s$ is continuous, non-decreasing and such that 
 \begin{equation}
 \label{approxtpsloc}
 \forall \, t,  a \geq 0 , \quad  \lim_{\varepsilon \rightarrow 0} \bbE \left[  \sup_{ 0\leq s \leq t} \left| \frac{1}{\varepsilon} \int_0^s dr 
\un_{\{ a < H_r \leq a+ \varepsilon \}} -L_s^a \right| \right] =0\; . 
\end{equation}
The process $(L^a_s,  s\geq 0)$ is called the {\it $a$-local time of $H$}. Recall that $I$ stands for the infinimum process of $X$. First, note that $L_{t}^{0}= -I_t$, $t \geq 0$. Next, observe that the support of the random Stieltjes measure $dL^{a}_{{\cdot }}$ is contained in the closed set $\{t \geq 0:H_t=a\}$. A general version of the {\it Ray-Knight theorem for $H$} asserts the following: For any $x\geq 0$, set $T_x=\inf \{ t\geq 0 \; :\;
X_t=-x\}$. Then, {\it the process $(L^a_{T_x} \; ;\; a \geq 0)$ is a distributed as a  CSBP with branching mechanism $\psi$ and initial state $x$} (see Le Gall and Le Jan \cite{LGLJ1} Theorem 4.2 and \cite{DuLG} Theorem 1.4.1). The CSBP $(L^a_{T_x} \; ;\; a\geq 0)$ admits a cadlag modification that is denoted in the same way to simplify notation. An easy argument deduced from the approximation (\ref{approxtpsloc}) entails that 
$\int_0^a L^b_{T_x}  \, db= \int_0^{T_x}\un_{\{ H_t \leq a\}} dt $. This remark combined with an elementary formula on CSBPs (whose proof can be found in Le Gall \cite{LG99}) entails that  
\begin{equation}
 \label{kappadef}
 \bbE \Big[ \exp \Big(  \! -\!\mu L^a_{T_x} -\lambda \!\! \int_0^{T_x} \!\!\! \un_{\{ H_t \leq a\} } dt  \, \Big) \Big]= \exp \! \big(\! - x \kappa_a (\lambda , \mu)\,  \big) \; ,  \quad a, \lambda , \mu \geq 0 ,
\end{equation}   
where $\kappa_a (\lambda , \mu)$ is the unique solution of the following differential equation 
\begin{equation}
\label{equakappa}
\kappa_0 (\lambda, \mu) = \mu  \quad {\rm and} \quad   \frac{\partial \kappa_a}{\partial a} (\lambda ,\mu ) = \lambda -\psi \big( \kappa_a (\lambda, \mu) \, \big)  \; , \quad a , \lambda , \mu \geq 0 .
\end{equation}

   It is possible to define the local times of $H$ under the excursion measure $N$ as follows. For any $b >0$, let us  set $v(b)=N( \sup_{^{t \in [0, \zeta ]}} H_t > b )$. Since $H$ is continuous, the Poisson decomposition (\ref{Poissheight}) implies that $v(b) < \infty$, for any $b >0$. 
It is moreover clear that $v$ is non-increasing and that $\lim_{\infty} v= 0 $. Then, for every $a\in (0, \infty)$, we define a continuous increasing process $(L^a_t,t\in [0, \zeta] )$, such that for every
$b \in (0,\infty)$ and for any $t\geq 0$, one has 
\begin{equation}
\label{localapprox}
\lim_{\varepsilon \rightarrow 0} \, 
N  \Big(  \un_{\{\sup H>b \}} \!\!  \sup_{ 0\leq s \leq t\wedge \zeta} \Big|  \frac{1}{\varepsilon} \int_0^s
dr 
\un_{\{ a-\varepsilon< H_r \leq a\}} -L_s^a \Big| \, \Big) =0.
\end{equation}
We refer to \cite{DuLG} Section 1.3 for more details. The process $(L^a_t,t\in [0, \zeta] )$ is the {\it $a$-local time of the excursion of the height process}. The Poisson decomposition (\ref{Poissheight}) then entails that 
\begin{equation}
\label{kappaexc}
N \left( 1- e^{ -\mu L^a_\zeta - \lambda \int_0^a \un_{\{ H_t \leq a \}} dt} \right) =
 \kappa_a (\lambda, \mu) \; , \quad a , \lambda , \mu , \geq 0 . 
\end{equation}
By taking $\lambda= 0$ in the previous display, we get $N ( 1- \exp (-\mu L^{_a}_{^\zeta} )\, ) = u(a, \mu)$, where $u$ is the solution of the integral equation  (\ref{integral}). This easily entails
\begin{equation}
\label{meanloc}
\forall a \geq 0 \; , \quad N (L^a_\zeta )= e^{-\alpha \, a} \; .
\end{equation}
Let us also recall from \cite{DuLG} the following formula  
\begin{equation}
\label{vvvequa}
\forall a>0 \; , \quad v(a)= N \big( \sup H_t \geq a \big) = N \left( L^a_\zeta \neq 0\right) \quad {\rm and} \quad \int_{v(a)}^\infty \frac{du}{\psi (u)} = a .    
 \end{equation}

\smallskip

\noi
{\bf L\'evy trees.}  We first define $\bR$-trees (or real trees) that are metric spaces generalising graph-trees. 
\begin{definition}
\label{errtreedef} Let $(T, \delta)$ be a metric space. It is a {\it real tree} iff the following holds true for any $\sigma_1, \sigma_1\in T$. 
\begin{description}
\item[{\bf (a)}] There is a unique isometry 
$f_{\sigma_1,\sigma_2}$ from $[0,\delta(\sigma_1,\sigma_2)]$ into $T$ such
that $f_{\sigma_1,\sigma_2}(0)=\sigma_1$ and $f_{\sigma_1,\sigma_2}(
\delta(\sigma_1,\sigma_2))=\sigma_2$. We denote by  $\lgeo \sigma_1,\sigma_2\rgeo$ the geodesic   joining $\sigma_1$ to $\sigma_2$. Namely, $\lgeo \sigma_1,\sigma_2\rgeo:=f_{\sigma_1,\sigma_2}([0,\delta (\sigma_1,\sigma_2)])$
\item[{\bf (b)}]   If $\inject$ is a continuous injective map from $[0,1]$ into
$T$, such that $\inject(0)=\sigma_1$ and $\inject(1)=\sigma_2$, then  we have
$\inject([0,1])=\lgeo \sigma_1,\sigma_2\rgeo$.
\end{description}
\noi
A rooted $\bR$-tree is a $\bR$-tree $(T,\delta)$ with a distinguished point ${\rm r}$ called the root.  \qed 
\end{definition}

 Among metric spaces, $\bR$-trees are characterized by the so-called {\it four points inequality}: $(T, \delta)$ is a $\bR$-tree iff it is connected and for any $\sigma_1,  \sigma_2,  \sigma_3,  \sigma_4  \in T$, 
\begin{equation}
\label{fourpoint}
\delta(\sigma_1, \sigma_2) + \delta(\sigma_3, \sigma_4) \leq \big(\delta(\sigma_1, \sigma_3) + \delta(\sigma_2, \sigma_4)\big) \vee  \big( \delta(\sigma_1, \sigma_4) + \delta(\sigma_2, \sigma_3)  \big) . 
\end{equation}
We refer to Evans \cite{EvStF} or to Dress, Moulton and Terhalle \cite{DMT96} for a detailed account on this property. The set of all compact rooted $\bR$-trees can be equipped with the {\it pointed Gromov-Hausdorff distance} that is defined as follows.   
\begin{definition}
\label{Hausdorffdist} {\bf (a)} Let $(E, \Delta)$ be a metric space. For any $x \in E$ and any subset $A \subset E $, we set $\Delta (x, A) = \inf_{y \in A} \Delta (x,y)$. Note that $\Delta ( \cdot , A)= \Delta ( \cdot , \bar{A})$ and that $\Delta (\cdot , A)$ is $1$-Lipschitz. For any $\varepsilon >0$, we set $A^{(\varepsilon)}=\{ x \in E: \Delta (x, A) \leq \varepsilon \}$ that is a closed subset of $E$. Then for any compact sets $K_1$, $K_2$ of $E$, we set 
$$\Delta_{{\bf H}} (K_1, K_2)= \inf \{ \varepsilon \in (0, \infty) : K_1 \subset K^{(\varepsilon)}_2 \; {\rm and} \; K_2 \subset K^{(\varepsilon)}_1 \} \; .$$ 
$\Delta_{{\bf H}}$ is a distance on the compacts sets of $E$ and we recall Blaschke's Theorem that asserts that the set of compact subsets of $E$ equipped with $\Delta_{{\bf H}} $ is a compact metric space when $(E, \Delta)$ is compact. 

\smallskip

\noi
{\bf (b)} Let $(T_{1}, \delta_1, {\rm r}_1)$ and $(T_{2}, \delta_2, {\rm r}_1)$ be two compact pointed metric spaces. The {\it pointed Gromov-Hausdorff distance} is then given by 
$$d_{{\bf GH}}(T_1,T_2)= \inf  \;  \; \Delta_{{\bf H}}\big( \, \inject_1(T_1),\inject_2 (T_2) \, \big) \vee \Delta \big( \, \inject_1 ({\rm r}_1), \inject_2 ({\rm  r}_2) \, \big) \; . $$  
where the infimum is taken over all the $(\inject_1,\inject_2, (E, \Delta))$, where $(E, \Delta) $ is a metric space and where $\inject_1 : T_1 \rightarrow E$ and $\inject_2 : T_2 \rightarrow E$ are isometrical embeddings. \cq 
\end{definition}
Obviously $d_{{\bf GH}}(T_1, T_2)$ only depends on
the root-preserving isometry classes of $T_1$ and $T_2$. In \cite{Gro}, Gromov proves that $d_{{\bf GH}}$ is a metric on the set of the equivalence classes of pointed compact metric spaces that makes it complete and separable. 
Let us denote by $\bbT$, the set of all 
equivalence classes of rooted compact real-trees. Evans, Pitman and Winter \cite{EvPitWin} prove that $\bbT$ is $d_{{\bf GH}}$-closed. Therefore, {\it $(\bbT, d_{{\bf GH}})$ is a complete separable metric space} (see Theorem 2 \cite{EvPitWin}). 

\smallskip

 Let us briefly recall how $\bR$-trees can be obtained via continuous functions. Recall from Notation \ref{excuheightdef} that $\ccC$ stands for the set of the continuous functions from $[0, \infty)$ to $\bR$ with compact support. Let $h \in \ccC$. To avoid trivialities, we also assume that $h$ is not constant to zero. Then, for every $s,t\geq 0$, we set
\begin{equation}
\label{pseudometric}
b_h(s,t)=\inf_{r\in[s\wedge t,s\vee t]}h(r) \quad {\rm and} \quad d_h(s,t)=h(s)+h(t)-2b_h(s,t).
\end{equation}
Clearly $d_h(s,t)=d_h(t,s)$. It is easy to check that $d_h$ satisfies the four points inequality, which implies that $d_h$ is a pseudo-metric.  We then introduce the equivalence relation
$s\sim_h t$ iff $d_h(s,t)=0$ (or equivalently iff $h(s)=h(t)=b_h(s,t)$) and we denote by $T_h$ the quotient set $[0,\zeta (h)]/ \sim_h$, where we recall that $\zeta (h)$ stands for the lifetime of $h$.  Standard arguments imply that $d_h$ induces a metric on $T_h$
that is also denoted by $d_h$ to simplify notation. We denote by
$p_h:[0,\zeta (h)]\rightarrow T_h$ the canonical projection. Since $h$ is continuous, so is $p_h$. This implies that $(T_h, d_h)$ is a compact and connected metric space that satisfies the four points inequality. It is therefore a compact $\bR$-tree. We then define the root $\rho_h$ of $(T_h, d_h)$ by $\rho_h= p_h (0)$. We shall refer to the rooted compact $\bR$-tree $(T_h, d_h, \rho_h)$ as to the {\it 
$\bR$-tree coded by $h$}.

  It shall be sometimes convenient to extend the canonical projection: we define $\bar{p}_h : [0, \infty) \rightarrow T_h$ by setting $\bar{p}_h (t)= p_h (t \wedge \zeta (h) ) $, $t \in [0, \infty)$. We next introduce the {\it mass measure on $T_h$}: We denote by $\ell$ the Lebesgue measure on $[0, \infty)$ and we denote by $\bm_h$ the measure on the Borel sets of $(T_h, d_h)$ induced by the measure $\ell$ restricted to $[0, \zeta (h)]$ via $p_h$. Namely, for any Borel subset $B$ of $T_h$, 
\begin{equation}
\label{massmeadeter}
 \bm_h (B)= \ell \big( p_h^{_{\; -1}} (B) \big)= \ell \big( [0, \zeta (h)] \cap  \bar{p}_h^{ _{ \; -1}} (B) \big) \; .
\end{equation}

  We next define {\it the $\psi$-L\'evy tree} as the tree coded by the $\psi$-height process $(H_t , 0 \leq t \leq \zeta)$ under the excursion measure $N$. To simplify notation, we set 
$$ (T_H, d_H, \rho_H , \bm_H)= (\cT, d , \rho, \bm) \; .$$
We also set $p= p_H: [0, \zeta] \rightarrow \cT$. Note that $\rho = p(0)$. Since $H_\zeta = 0$ and since $H_t >0$, for any $t \in (0, \zeta)$, $\zeta$ is the only time $t$ distinct from $0$ such that $p(t)= \rho$.

  A point  $\sigma \in \cT$ is called a {\it leaf} if it is distinct from the root and if the open set $\cT \backslash \{ \sigma \}$ is connected. We denote by ${\bf Lf  } (\cT)$ the set of leaves of $\cT$. We also define the {\it skeleton of $\cT$} by ${\bf Sk} (\cT)= \cT \backslash{\bf Lf  } (\cT) $. One can show that 
 \begin{equation}
\label{massskel}
\textrm{$N$-a.e.$\quad \overline{{\bf Sk}} (\cT)= \cT$, $\, \bm$ is diffuse and $\; \bm \big( {\bf Sk} (\cT) \big) = 0 $.}
\end{equation}
This easily implies the following characterisation of leaves in terms of the height process: $N$-a.e.$\;$for any $t\in (0, \zeta)$, 
\begin{equation}   
\label{Hleafcarac}
 p(t) \in {\bf Lf} (\cT)  \quad  \Longleftrightarrow \quad  \forall \varepsilon  > 0 \, , \;  \inf_{s \in [t-\varepsilon , t]} \!\! H_s < H_t  \quad  {\rm and} \;  \inf_{s \in [ t,  t+\varepsilon ]} \!\! H_s  \;  < \;  H_t \;  . 
\end{equation}
For any $a \in (0, \infty)$, the {\it $a$-local time measure} $\ell^a$ is the measure induced by $dL^a_{\cdot}$ via $p$. Namely, 
$$ \langle \ell^a , f \rangle = \int_0^\zeta \!\! \! dL^a_s \,  f(p (s)) \; , $$
for any positive measurable application $f$ on $\cT$. Let us mention that the topological support of $\ell^a$ is included in  the $a$-level set $\cT (a)= \{ \sigma \in \cT: d(\rho ,\sigma)= a \}$ and note that the total mass $\langle \ell^a\rangle$ of $\ell^a$ is equal to $L^{_a}_{^\zeta}$. Moreover, observe that $\cT(a)$ is not empty iff $\sup H \geq a$. Then, (\ref{vvvequa}) can be rewritten as follows. 
\begin{equation}
\label{treevvvequa}
\forall \, a>0 , \quad v(a)= N \big( \cT(a) \neq \emptyset  \big) = N \left( \ell^a  \neq 0\right)  \; .    
 \end{equation}
As already mentioned, the $a$-local time measure $\ell^a$ can be defined in a purely metric way by (\ref{geolocapprox}) and there exists a modification of 
$a \mapsto \ell^a$ that is $N$-a.e.$\;$cadlag for the weak topology on the space of finite measures on $\cT$. 
  
\smallskip

For any $h \in \ccC$, denote by $\bar{T}_h$ the root-preserving isometry class of $(T_h, d_h, \rho_h)$ that belongs to $\bbT$. Lemma 2.3 \cite{DuLG2} asserts that $h\in \ccC \mapsto \bar{T}_h \in \bbT$ is Borel-measurable.  We then define $\Theta_\psi$ as the "distribution" of $\bar{\cT}$ when $\cT$ is under $N$. We have stated the main results of the paper under $\Theta_\psi$ because it is more natural and because $\Theta_\psi $ has an intrinsic characterization as shown by Weill \cite{Weill}. However, each time we make explicit computations with L\'evy trees, we have to work with random isometry classes of compact $\bR$-trees, which causes technical problems (mostly measurability problems). To avoid these unnecessary complications during the intermediate steps of the proofs, {\it we prefer to work with the specific compact rooted $\bR$-tree $(\cT, d, \rho)$ coded by the $\psi$-height process $H$ under $N$ rather than to directly work under $\Theta_\psi$}.

\smallskip

\noi
{\bf The branching property.} We now describe the distribution of the subtrees above level $b$ in the L\'evy tree. More precisely, we consider the excursions above level $b$ of the height process $H$ under $N$. Let us fix $b\in (0, \infty)$. We denote by $(g^{b}_{j}, d^{b}_{j})$, $j\in \cI_b$, 
the connected components of the open set $\{t \geq 0:H_t>b\}$. For any $j\in \cI_b$,
we denote by $H^{{b, j}}$ the corresponding excursion of $H$ defined by $H^{{b, j}}_s=H_{(g^b_j+s)\wedge d^b_j}-b$, $ s\geq 0$. This has to be interpreted in terms of the tree as follows. Recall that $\bar{B} (\rho , b)$ stands for the closed ball in $\cT$ with center $\rho$ and radius $b$. Observe that the connected components of the open set $\cT \backslash \bar{B} (\rho , b)$ are the subtrees $\cT_{j}^{{b, o}}:= p((g^{b}_{j}, d^{b}_{j}))$, $j \in \cJ_b$. The closure $\cT_{j}^{b}$ of $\cT_{j}^{{b, o}}$ is simply $\{ \sigma^{b}_{j} \} \cup \cT_{j}^{{b, o}}$, where $\sigma^{b}_{j} = p(g^{b}_{j})= p(d^{b}_{j})$ is the points on the $b$-level set $\cT(b)$ at which $\cT^{{b, o}}_{j}$ is grafted. Observe that the rooted compact $\bR$-tree $(\cT^{b}_{j}, d, \sigma^{b}_{j})$ is isometric to the tree coded by $H^{{b, j}}$.

  We then define $\tilde{H}^b_s=H(\zeta \wedge \tau^{b}_{s})$, where $s \mapsto \tau^{b}_{s}$ is given by  
$$\forall s \geq 0 \, , \quad \tau^{b}_{s}= \inf \Big\{  t\geq 0: \, \int_0^t \!\! dr\,\un_{\{H_r\leq b\}}>s \Big\}\, ,$$  with the usual convention $\inf \emptyset = \infty$. 
The process $\tilde{H}^b$ is the {\it height process below $b$} and the rooted compact $\bR$-tree 
$(\bar{B} (\rho , b), d, \rho)$ is isometric to the tree coded by $\tilde{H}^b$. We denote by $\cG_b$ the 
sigma-field generated by $\tilde{H}^b $ augmented by the $N$-negligible sets.We see from (\ref{localapprox}) that $L^{b}_{\zeta}$ is measurable with respect
to $\cG_b$. We next define the probability measure $N_b$ on $\cC^0$ by  
\begin{equation}
\label{condiprob}
N_b = N (\, \,  \cdot \,\,  | \, \sup H>b)
\end{equation}
and we introduce the following point measure on $[0, \infty) \times \cC^0$:
\begin{equation}
\label{branchprop}
\cM_b =  \sum_{^{j\in \cI_b}} \delta_{\big( L^b_{^{g^{_b}_{^j}}} \; , \, H^{b,j}  \big)}
\end{equation}
{\it The branching property at level $b$ then asserts that under $N_b$, conditionally given $\cG_b$, $\cM_b $ is distributed as a Poisson point measure with intensity} 
$\un_{[0,L^{b}_{\zeta} ]}(x)dx\otimes  N(dH)$ (see \cite{DuLG} Proposition 1.3.1).  Let us mention that it is possible to rewrite intrinsically the branching property under $\Theta_\psi$: see \cite{DuLG2} Theorem 4.2, for more details. As already mentioned Weill \cite{Weill} shows that the branching property characterizes L\'evy trees.

\smallskip

\noi
{\bf Spinal decomposition.} 
Let us introduce an auxiliary probability space $(\Omega, \cF , \bP)$ that is rich enough to 
carry the various independent random variables we shall need. Let $Y= ( W_{t}, V_{t} )_{t \geq 0}$ be a bivariate subordinator on $(\Omega, \cF, \bP)$ with initial value $Y_0= (0,0)$. Namely, $Y$ is a cadlag process with independent and homogeneous nonnegative increments. Its distribution is characterised by its Laplace exponent given by 
$$ -\frac{1}{t}  \log \bE \big[ \exp (  -\lambda  W_{t}    -\mu V_{t}   ) \big]   =  \frac{{\psi^* (\lambda)-\psi^*(\mu)}}{{\lambda-\mu  }}     \; , $$
where $\psi^* (\lambda)= \psi (\lambda)-\alpha \lambda$. If $\lambda= \mu$, the right member has to be interpreted as the derivative $(\psi^*)^\prime (\lambda)$. Denote by $\overline{W}$ and $\overline{V}$ the right-continuous inverses of $W$ and $V$: 
$$ \overline{W} (r) = \inf \{ t \in [0, \infty); W_t >r \} \quad {\rm and} \quad  \overline{V}(r) = \inf \{ t \in [0, \infty); V_t >r \} \; .$$
Note that $W$ and $V$ are two subordinators with Laplace exponent $\psi^* (\lambda)/\lambda$. Since (\ref{extinction}) implies that $\beta >0$ or $\int_{(0, 1)} r\pi (dr) = \infty$,  $\bP$-a.s.$\, W$ and $V$ are increasing. 
Thus $\overline{W}$ and $\overline{V}$ are $\bP$-a.s.$\;$continuous. 

Let $(X^{_{(1)}}_{t} )_{ t \geq 0}$ and $(X^{_{(2)}}_t )_{t \geq 0}$ be two independent real valued L\'evy processes defined on $(\Omega, \cF , \bP)$ whose common distribution is $\bbP$. Thus, their initial value is $0$, their Laplace exponent is $\psi$. 
 {\it We moreover assume that $(X^{(1)}, X^{(2)} )$ is independent of $Y$}. We denote by $H^{(1)}$ and $H^{(2)}$ the height processes obtained respectively from $X^{(1)}$ and from $X^{(2)}$. Thus, $H^{(1)}$ and $H^{(2)}$ are two independent $\psi$-height processes and $(H^{(1)}, H^{(2)})$ is independent from $Y$. From $(H^{(1)}, H^{(2)})$ and $Y$, we derive two processes as follows. For any $t \in [0, \infty)$, we set 
\begin{equation}
\label{starundeudef}
H^{{*(1)}}_{t} = H^{{(1)}}_{t} -\overline{W} (-I^{{(1)}}_{t} ) \quad {\rm and } \quad H^{{*(2)}}_{t} =H^{{(2)}}_{t} -\overline{V} (-I^{{(2)}}_{t} )\; , 
\end{equation}
where $I^{_{(1)}}_{^t}= \inf_{^{s \in [0, t]}} X^{_{(1)}}_{^s} $ and $I^{_{(2)}}_{^t}= \inf_{^{s \in [0, t]}} X^{_{(2)}}_{^s} $. Observe that $H^{*(1)}$ and  $H^{*(2)}$ are continuous, possibly negative, that $H^{{*(1)}}_{0}=H^{{*(2)}}_{0}= 0$ and that 
\begin{equation}
\label{infistar}
\forall \, t \in [0, \infty) \, , \quad \inf_{{s \in [0, t]}} H^{{*(1)}}_{s}=  -\overline{W}(-I^{{(1)}}_{t} ) \quad {\rm and} \quad \inf_{{s \in [0, t]}} H^{{*(2)}}_{s}=  -\overline{V}(-I^{{(2)}}_{t} ) \; .
\end{equation}
Next, for any $a \in [0, \infty)$, we set 
$$T^{{(1)}}_a = \inf \{ t \in [0, \infty) \, ; \;  H^{{*(1)}}_{ t} = -a \} \quad {\rm and} \quad T^{{(2)}}_a = \inf \{ t \in [0, \infty) \, ; \; H^{{*(2)}}_{ t} = -a \} . $$
Note that $T^{_{(1)}}_{a}= \inf \{ t \geq 0 \, ; \;  X^{_{(1)}}_{^t} = -W_a \}$ and $T^{_{(2)}}_{a}= \inf \{ t \geq 0 \, ;\;  X^{_{(2)}}_{^t} = -V_a \}$. We next set 
\begin{equation}
\label{aundeudef}
 H^{(a, 1)} = \big( a+ H^{_{*(1)}}_{ t \wedge T^{_{(1)}}_a } \, , \, t \geq 0 \big)  \quad {\rm and} \quad H^{(a, 2)} = \big( a+ H^{_{*(2)}}_{ t \wedge T^{_{(2)}}_a } \, , \, t \geq 0 \big). 
\end{equation}
They are nonnegative continuous processes  with compact support and with respective lifetimes $T^{_{(1)}}_a$ and $T^{_{(2)}}_a$
Let us now consider the height process  $H= (H_t, t \geq 0)$ under $N$. 
For any $t \geq 0$, we set 
\begin{equation}
\label{hatcheck}
\hat{H}^{t}:= ( H_{(t- s)_+}, s \geq 0)  \quad {\rm and} \quad \check{H}^{t}:= (H_{t+s} , s \geq 0) \; , 
\end{equation}
where, $(\, \cdot )_+$ stands for the positive part function.
Then, for any bounded measurable function $F:\cC^0 \times \cC^0  \rightarrow [0, \infty)$, one has 
\begin{equation}   
\label{BismutHeight}
N \Big( \!\!  \int_{0}^\zeta  \!\!\!  F \big( \hat{H}^{t} , \check{H}^{t} \big)  \,  dt  \Big) 
= \int_0^\infty \!\!\! e^{-\alpha a} \bE  \Big[ \, F\big(H^{(a, 1)} , H^{(a, 2)}  \big) \,  \Big] \, da \; .  
\end{equation}
In the Brownian case, this decomposition is equivalent to Bismut decomposition. As already mentioned, 
this decomposition is a consequence of Lemma 3.4 \cite{DuLG2} (see also \cite{DuLG} Chapter 1 or  Lemma 3.2 \cite{Du4}, and see \cite{DuLG4} for further applications).

   We first use (\ref{BismutHeight}) to prove the following zero-one law that is needed in the proof of Lemma \ref{technimesu} in Appendix \ref{appen}. For any $\eta \in (0, \infty)$, we define $R_\eta: \cC^0 \rightarrow \cC^0$ by setting $R_\eta H= (H_{s\wedge \eta }-H_\eta, s \geq 0)$. Note that $R_\eta$ is continuous. 
We next fix a sequence $\eta_n \in (0, \infty)$, $n \geq 0$,  that decreases to $0$. We also fix a sequence of Borel-measurable functions $G_n: \cC^0 \rightarrow [0, \infty]$, $n \geq 0$. We then set $ G (H)= \liminf_{n \rightarrow \infty} G_n ( R_{\eta_n} H ) \in [0, \infty] $ that is measurable from $\cC^0$ to $[0, \infty]$. 
\begin{lemma}
\label{zerooneprelim} There exists a constant $C\in [0, \infty]$ such that 
$ N \big( \int_0^\zeta \! \un_{\{ G ( \check{H}^t  ) \neq C  \}}  dt \big)  \!= \!0 $. 
 \end{lemma}
\noi
{\bf Proof:} (\ref{BismutHeight}) implies $N \big( \int_0^\zeta \! \un_{\{ G ( \check{H}^t  ) \neq C  \}}  dt \big)= \int_0^\infty e^{-\alpha a} \bP ( G( H^{(a, 2)})\neq C)da $, for any $C \in [0, \infty]$. 
For any $\eta, a \in (0, \infty)$, we set $\tau (\eta, a)= \eta\wedge T^{_{(2)}}_{a}$. Then observe that for any $s \geq 0$, we have $ R_\eta  H^{(a, 2)}  (s) =   H^{_{(2)}}_{^{s\wedge \tau (\eta, a)}} \!- \!H^{_{(2)}}_{^{\tau (\eta, a)}} \! - \! \overline{V} (-I^{_{(2)}}_{^{s \wedge \tau (\eta, a)}} \,  )\!+ \! \overline{V} (-I^{_{(2)}}_{^{ \tau (\eta, a)}} ) $.
Note that $G(H^{(a, 2)})$ is $\bP$-a.s.$\;$equal to a random variable that does not depend on $a$ and that is measurable with respect to the tail sigma-field at $0+$ of the Feller process 
$ (X^{{(1)}}, X^{{(2)}}, W , V)$. By Blumenthal's zero-one law, there exists a constant $C\in [0, \infty]$ such that for any $a$, $ \bP ( G( H^{(a, 2)})\neq C)= 0$, which implies the desired result. \cqfd

 \smallskip

We now recall from \cite{DuLG2} a Poisson decomposition of $H$ that is a consequence of (\ref{BismutHeight}). For any continuous function $h\in \cC^0$, we define the point measure $\cN (h)$ as follows. Set $\underline{h} (t)= \inf_{[0, t]} h$ and denote by $(g_i, d_i)$, $i\in \cI(h)$ the excursion intervals of $h-\underline{h} $ away from $0$ that are the connected components of the open set $\{ t \geq 0: h(t)-\underline{h} (t) >0 \}$. For any $i \in \cI (h)$, set 
$h^i (s)  = ((h- \underline{h}) ( (g_i +s) \wedge d_i) \, , \, s \geq 0)$ and define 
$$ \cN (h)= \sum_{ i\in \cI (h)} \delta_{\big( h(0)-h(g_i)\,  , \,  h^i \big) } $$
that is a point measure on $[0, \infty) \times \cC^0$. Recall that $H= (H_t)_{t \geq 0}$ stands for the excursion of the height process $H$. For any $t \in [0, \zeta]$, we set  
\begin{equation}
\label{spinaldef}
\cN_t =\cN( \hat{H}^{t})+ \cN ( \check{H}^{t}) := \sum_{j \in \cJ_t } \delta_{(r^t_j, H^{\bullet \, t,j} )} \; \, . 
\end{equation}
Recall from (\ref{starundeudef}) the definition of $H^{*(1)}$ and $H^{*(2)}$. Then, we also set, 
\begin{equation}
\label{Nstaradef}
\cN^*= \cN ( H^{*( 1)}) + \cN (H^{*(2)} ):=  \sum_{j \in \cI^*}  \, \delta_{ (  r^*_j , \, H^{*j} )  }
\end{equation}
By definition of $H^{(a, 1)}$ and $H^{(a,2)}$,  it is easy to check that 
$$ \cN_a^*:=  \cN ( H^{(a, 1)}) + \cN (H^{(a, 2)} ):=  \sum_{j \in \cI^*} \un_{[0, a]} (r^*_j) \, \delta_{ (  r^*_j , \, H^{*j} )  } . $$
Then, (\ref{BismutHeight}) implies that 
\begin{equation}      
\label{ancdecomp}
N \Big( \int_{0}^\zeta \!\!\!   \; F \big(  \cN_t  \big) \, dt  \Big) =  \int_0^\infty \!\!\! \! \!  \,  e^{-\alpha a } \,   \bE \left[ F (\cN^*_a ) \right] \, da \; . 
\end{equation}
We shall refer to this identity as to the {\it spinal decomposition of $H$ at a random time} (see \cite{DuLG2} Lemma 3.4). Let us briefly explain the distribution of $\cN^*$ under  $\bP$: Recall that for any $k \in \{ 1, 2\}$, $H^{(k)}$ is the height process associated with $X^{(k)}$. We denote by $H^{(k),j}$, $j \in \cJ_k$, the 
excursions of $H^{(k)}$ above $0$, and we denote by $(g(k,j), d(k,j))$ the corresponding excursion intervals.  As a consequence of (\ref{infistar}) and of the definition of $H^{*(1)}$, the atoms of $\cN (H^{*( 1)})$ are the points
$( \overline{W} (-I^{_{(1)}}_{^{g(1,j)}}) \,  , \,  H^{(1),j}) $, $j \in \cJ_1$. Similarly, the atoms of $\cN (H^{*(2)})$ are the points
$( \overline{V} (-I^{_{(2)}}_{^{g(2,j)}}) \,  , \,  H^{(2),j}) $, $j \in \cJ_2$. We then set $U_r= W_r+V_r$, $r \in [0, \infty)$. Then $U= (U_r)_{r \geq 0}$ is a subordinator with Laplace exponent $(\psi^*)^\prime$ and it is easy to check that $dU_r= dW_r + dV_r$. Since the measure induced by the Lebesgue measure on $[0, \infty)$ via $\overline{W}$ 
(resp.$\;$via  $\overline{V}$) is the random Stieltjes measure $ d W_r $ (resp.$\, d V_r $), (\ref{Poissheight}) implies that for any measurable function 
$\Phi : [0, \infty) \times \cC^0 \rightarrow [0, \infty]$, 
\begin{equation}
\label{conditNstar}
\bE \big[ \exp \big( \!- \! \langle \cN^* , \Phi \rangle  \big) \, \big| \, U \big] = \exp\Big( \!- \! \int_{0}^\infty \!\!\! dU_r \,   N \big( 1-e^{\Phi(r, H)}\big) \, \Big) \; .
\end{equation}
Thus, {\it the law of $\cN^*$ conditionally given $U$ is that of a Poisson point measure with intensity $dU_r \otimes N(dH)$}.

Let us briefly interpret this decomposition in terms of the $\psi$-L\'evy tree $\cT$ coded by $H$ under $N$. Choose $t \in (0, \zeta)$ and set $\sigma = p(t) \in \cT$.  Then the geodesic $\lgeo \rho, \sigma \rgeo$ is interpreted as the ancestral line of $\sigma$. Let us denote by $\cT_{j}^{o}$, $j \in \cJ$, the connected components of the open set $\cT   \backslash \lgeo \rho , \sigma \rgeo$ and denote by $\cT_{{j}}$ the closure of 
$\cT_{j}^{o}$. Then, there exists a point $\sigma_j \in \lgeo \rho , \sigma \rgeo$ such that $\cT_{{j}} = \{ \sigma_{{j}} \} \cup \cT_{j}^{o}$. Recall that $(r^{t}_{j} , H^{{{\bullet \, t,j}}})$, $j \in \cJ_t$ are the atoms of $\cN_t$ as defined by (\ref{spinaldef}). Then, for any $j \in \cJ$, there exists a unique $j^\prime \in \cJ_{t}$ such that $ d(\sigma ,\sigma_{{j}})= r^{t}_{{j^\prime}} $ and such that the rooted compact $\bR$-tree $(\cT_{{j}}, d, \sigma_{{j}})$ is isometric to the tree coded by $H^{{{\bullet \, t,j^\prime}}}$.

   We now apply (\ref{ancdecomp}) to compute the mass measure of balls whose center is chosen according to the mass measure $\bm$. Let $t \in (0, \zeta)$. We first compute $\bm (\bar{B}(p(t), r))$ in terms of  $\cN_t$ as follows. Recall notation $b(s,t)= \min_{[s\wedge t, s\vee t]} H$. By (\ref{Hleafcarac}), $N$-a.e.$\;$ for all $s, t\in (0, \zeta)$ such that $H_s=b(s,t)$ with $s \neq t$, we have $p(s) \in {\bf Sk} (\cT)$. 
By (\ref{massskel}), $N$-a.e.$\, \bm ({\bf Sk} (\cT))= 0$. Consequently, $N$-a.e.$\;$for every $r\in (0, \infty)$ and every 
$t \in (0, \zeta)$, we have 
\begin{eqnarray}
\label{massballfixt}
\bm \big( \bar{B} (p(t) , r  ) \big) &= & \int_{0}^\zeta \un_{\{ d(s,t) \leq  r \}} ds = \int_0^\zeta 
\un_{\{ 0 < H_s  - b(s,t ) \leq r -  H_t+b(s,t)     \}} \nonumber \\
& =& \sum_{j \in \cJ_t} \un_{[ 0 \, , \, r]} (r^t_j)  \int_0^{\zeta^t_j} \un_{\{ H^{\bullet \, t,j}_s \leq r-r^t_j\}} \, ds , 
\end{eqnarray}
where $\zeta^t_j$ stands for the lifetime of the path $H^{\bullet \,  t,j}$. For any $a \in (0, \infty)$ and for any $r \in [0,\infty)$, we next set 
\begin{equation}
\label{Mstardef}
M^*_r (a)= \sum_{j \in \cI^*} \un_{[0\,,\,  r \wedge a]} (r^*_j)  \int_{0}^{\zeta^*_j} \un_{\{ H^{*j}_s \leq r-r^*_j  \}} ds  \; , 
 \end{equation}    
where $\zeta^*_j$ stands for the lifetime of the path $H^{*j}$. Then, $(M^*_r (a) , r \geq 0) $ is a cadlag increasing process defined on $(\Omega, \cF, \bP)$. The spinal decomposition (\ref{ancdecomp})  entails the following key formula that is used in the proof of Theorem \ref{densityth}: For any bounded measurable $F: \bbD ([0, \infty) , \bR) \rightarrow [0, \infty )$, we have 
\begin{eqnarray}
\label{keymass}
N \left( \int_{\cT} \!\!    F  \! \left(   \, \bm \big( \bar{B}(\sigma , r) \, \big)  , \, r \geq 0 \,  \right)  \bm ( \!\!\; d\sigma \! \! \; )\! \right)  & =&  N \left( \int_0^\zeta \!\!\!  \, F \! \left( \,   \bm \big( \bar{B}(p(t) , \! r) \, \big)  , \, r \geq 0 \,  \right)   dt \right) \nonumber \\
& =&  \int_0^\infty \!\!\! \!\!  e^{-\alpha a } \,  \bE \left[  F \big( M^*_r  (a) , \, r \geq 0  \, \big)\right] \, da \; .
\end{eqnarray}

\subsection{Exponents.}
\label{exponentsec}
In this section we relate several power exponents associated with $\psi$ to properties of the gauge function $g$ that is derived from $\psi$ by (\ref{gaugedef}). Let us start with some notation. Let $\phi: [0, \infty) \rightarrow [0, \infty)$ be a continuous increasing function. We agree on the following conventions: $\sup \emptyset = 0$ and $\inf \emptyset = \infty $, and we define the following exponents that compare $\phi$ with power functions at infinity. 
\begin{description}
\item{(a)}  $\gamma_{\phi} := \sup \{ c \geq 0 :  \lim_{\lambda \rightarrow \infty} \phi (\lambda) \lambda^{-c} = \infty  \} $ is the {\it lower exponent of $\phi$ at $\infty$}. 
\item{(b)} $\eta_{\phi}: = \inf \{ c \geq 0 : \lim_{\lambda \rightarrow \infty} \phi (\lambda) \lambda^{-c} = 0 \} $  is the {\it upper exponent of $\phi$ at $\infty$}. 
\item{(c)} $ \delta_{\phi}: = \sup \, \{ c \geq 0  :  \exists \, C \!\in \!(0, \infty)\; \textrm{such that} \; C \phi (\mu) \mu^{-c}  \leq \phi (\lambda) \lambda^{-c}  \, , \, 1 \leq \mu \leq \lambda  \} $. 
\end{description}

\noi
{\bf Exponents for subordinators.} Let us assume that $\phi$ is the Laplace exponent of a subordinator with drift ${\rm d}$ and L\'evy measure $\nu$:
$$ \phi (\lambda ) = {\rm d} \, \lambda + \int_{(0, \infty)} \!\! (1-e^{-\lambda r}) \, \nu (dr ) \; , \quad \lambda \geq 0 . $$
Then, we have $0 \leq \delta_\phi \leq  \gamma_\phi \leq \eta_\phi\leq  1$. Recall that $\lim_{\lambda \rightarrow \infty} \phi (\lambda )/ \lambda = {\rm d}$. If ${\rm d}>0$, then 
$\delta_\phi= \gamma_\phi= \eta_\phi= 1$. If ${\rm d}=0$, the exponents can be expressed in terms of the 
L\'evy measure as follows: For any $x \in [0, \infty)$,  we set
\begin{equation}
\label{Jdef}
J_\phi (x) = \int_0^x \nu \big( (u, \infty) \big) \, du= \int_{(0, x]} r \, \nu (dr) + x \nu \big( (x, \infty) \big) \; .
\end{equation}
By standard results (see Bertoin \cite{Be} Chapter III), if ${\rm d}= 0$, there exist two universal constants $k_1, k_2 \in (0, \infty)$ such that 
\begin{equation}
\label{phiJ}
k_1 \, \lambda J_\phi(1/ \lambda) \leq \phi (\lambda)  \leq k_2 \,  \lambda J_\phi(1/ \lambda ) \; , 
\quad \lambda > 0  \; . 
\end{equation}
This easily implies that $\gamma_\phi = \sup \{ c \geq 0 : \, \lim_{0+} x^{c-1} J_\phi (x) = \infty  \}$, that $\eta_\phi = \inf \{ c \geq 0 : \, \lim_{0+} x^{c-1} J_\phi (x) = 0  \} $ and that 
\begin{equation}
\label{Jdeltaexpon}
\delta_\phi = \sup  \{ c \geq 0 : \, \exists \, C \!\in \! (0, \infty)\; \textrm{s.t.} \; C y^{c-1} J_\phi (y)
\leq x^{c-1} J_\phi (x) , \,  0 < x\leq y \leq 1 \} . 
\end{equation}

\noi
{\bf Exponents for $\psi$.} Let $\psi$ be of the form (\ref{LevyKhin}). Set $\widetilde{\psi} (\lambda) = \psi (\lambda ) / \lambda$. It is easy to show that for any $\lambda  \in [ 0, \infty)$, 
$$ \psi' (\lambda)= \alpha \!+ \!2 \beta \lambda \!+ \!\! \int_{(0, \infty)} \!\!\!  \!\!\! \!\!\!(1-e^{\lambda r})\,  r\pi (dr) \quad 
{\rm and} \quad \widetilde{\psi} (\lambda)= \alpha \!+ \! \beta \lambda \!+\!\!  \int_{(0, \infty)}  \!\!\! \!\!\!\!\!\! (1-e^{\lambda r})\,  \pi \big([r, \infty) \big)  \, dr . $$
Thus, $\psi^\prime $ and $\widetilde{\psi}$ are Laplace exponents of subordinators. Recall that the reciprocal $\psi^{-1} $ of $\psi$ is the Laplace exponent of a subordinator. 
Thus, $\varphi= \psi^\prime \circ \psi^{-1}$ is also the Laplace exponent of a subordinator. Note that 
$1/ \varphi$ is the derivative of $\psi^{-1}$. Note that $\psi$ is convex and that $\psi^\prime$, $ \widetilde{\psi}$, 
$\psi^{-1}$ and $\varphi$ are concave. In particular, this implies $\widetilde{\psi}(2\lambda) \leq 2 \widetilde{\psi} (\lambda)$ and the following 

\begin{equation}
\label{convexpsi}
\psi (2\lambda ) \leq 4 \psi (\lambda) \; , \quad     \widetilde{\psi} (\lambda)  \leq \psi'(\lambda) \leq 4 \widetilde{\psi} (\lambda)  \quad {\rm and }\quad   \;  \frac{\lambda}{\psi^{-1} (\lambda) }  \leq \varphi (\lambda) \leq \frac{4\lambda}{\psi^{-1} (\lambda) }. 
\end{equation}
To simplify notation we set 
$$\gamma := \gamma_{\psi} \, , \; \eta:= \eta_{\psi} \quad {\rm  and} \quad \delta:= \delta_\psi \; .$$ 
We clearly have $1 \leq \delta \leq  \gamma \leq \eta \leq 2$. In general, $\gamma$ and $\eta$ are distinct but they coincide if $\psi$ is regularly varying at $\infty$. As a direct consequence of (\ref{convexpsi}) we have $\delta_{\widetilde{\psi}} =\delta_{\psi'} = \delta -1$, $\gamma_{\widetilde{\psi}}=\gamma_{\psi'}=  \gamma-1$ and $\eta_{\widetilde{\psi}} =\eta_{\psi'}= \eta -1 $. Moreover, we get $\delta_\varphi= (\delta -1)/\delta$, $\gamma_{\varphi}=  (\gamma -1)/\gamma$ and $\eta_{\varphi}= (\eta -1)/\eta$.

  Recall the definition of the gauge function $g: (0, r_0) \rightarrow (0, \infty)$ that is derived from $\psi$ 
by (\ref{gaugedef}). The function $g$ is clearly continuous. For any $r \in (0, r_0)$, we set $a _r = \varphi^{-1} ( \frac{1}{r}  \log \log \frac{1}{r}  ) $. First, observe that $a_r $ increases to $\infty$ when $r$ decreases to $0$. Next, recall that since $\varphi$ is the Laplace exponent of a subordinator, $\lambda \mapsto 
\varphi( \lambda ) / \lambda$ is decreases. Thus, $ r \in (0, r_0) \mapsto g(r)= r  \varphi (a_r) / a_r$ is continuous, increasing and it goes to $0$ when $r$ goes to $0$. The following lemma relates the doubling condition (\ref{doubling}) for $g$ to the exponent $\delta$ of $\psi$. 
\begin{lemma}
\label{doublinglemma}
Assume that $\psi$ is of the form (\ref{LevyKhin}) and that it satisfies $(\ref{extinction})$. Then, the following assertions hold true. 
\begin{itemize}
\item[(i)] The gauge function $g$ satisfies the doubling condition (\ref{doubling}) iff $\delta >1$. 
\item[(ii)] If $\psi$ is regularly varying at $\infty$ with exponent $c>1$, then $\delta= \gamma= \eta=c$ and $g$ satisfies the doubling condition (\ref{doubling}). 
\end{itemize}
\end{lemma}
\noi {\bf Proof:}  We first assume that $\delta >1$. Then, $\delta_\varphi = (\delta -1)/\delta \in (0, 1)$. Let $c\in (0, \delta_\varphi)$. There exists $C\in (0, \infty)$ such that $ C \varphi (a) a^{-c}  \leq \varphi (b) b^{-c}$, for any $ 1 \leq a \leq b$. If we take $a= \varphi^{-1} (u)$ and $b= \varphi^{-1} (v)$, with $ \varphi (1) \leq u \leq  v$, then we get 
\begin{equation}
\label{auxidoublemma}
 \frac{u^{1/c}}{\varphi^{-1} (u)} \leq C^{-1/c} \cdot  \frac{v^{1/c}}{\varphi^{-1} (v)}\; , \quad \varphi (1) \leq u \leq v \; . 
 \end{equation}
Let $r_1 \in (0, r_0)$ be such that $ \varphi (1) \leq \frac{1}{2r}\log \log\frac{1}{2r} \leq \frac{1}{r}\log \log \frac{1}{r} $ for any  $r \in (0, r_1)$. Apply (\ref{auxidoublemma}) with $u=  \frac{1}{2r}\log \log\frac{1}{2r}$ and $v= \frac{1}{r}\log \log \frac{1}{r}$ to get 
$$ g (2r) \;  \leq \left( \frac{2}{C}\right)^{1/c} \left(\frac{ \log  \log \frac{1}{r} }{ \log \log \frac{1}{2r}} \right)^{\!\!\frac{1}{c} -1} \! \! \! \! g (r) \; , \quad r \in (0, r_1) \; , $$
which easily entails the doubling property (\ref{doubling}) for $g$ since $\frac{1}{c} -1 >0$.

  Conversely, let us assume that $g$ satisfies a doubling  property: there exists $C^\prime \in (1, \infty)$ such that 
$$ \frac{\log\log \frac{1}{2r}}{\varphi^{^{-1}} \! \! \left( \frac{1}{2r} \log \log \frac{1}{2r}\right)  } \leq C^\prime    \frac{\log\log \frac{1}{r} }{\varphi^{^{-1}} \! \! \left( \frac{1}{r} \log \log \frac{1}{r}\right)  }  \; , \quad r \in (0 , r_0/2)\; .$$ 
The previous inequality, combined with an easy argument, implies there exists $C>1$ and $u_0 >0$ such that $ \varphi^{-1} (u) \leq C \varphi^{-1} (u/2) $, for any  $u \geq u_0$. The previous inequality entails that $ 2\varphi (v) \leq \varphi (C v)$ for any $ v \geq v_0:= \max (1, \varphi^{-1} (u_0 /2) )$.  
We next set $c = \log(2)/\log(C) $ that  is strictly positive since $C >1$. For any $\lambda \geq v \geq v_0$, we denote by $n (v, \lambda)$, the integer part of $\log(\lambda / v) / \log(C)$. Namely, $ C^{n(v, \lambda)} v \leq \lambda < C^{n(v, \lambda)+1} v $. This implies 
 $$ \frac{_1}{^2} \cdot \lambda^c v^{-c} \,  \varphi (v)  \leq  2^{n(v, \lambda)} \varphi (v) \leq \varphi ( C^{n(v, \lambda)}  v) \leq \varphi ( \lambda) \; ,$$
which implies $\delta_\varphi >0$ and thus, $(1-\delta_\varphi)^{-1}= \delta >1$. This completes the proof of $(i)$. 

\medskip

 The second point of the lemma is a direct consequence of a theorem due to Matuszwska \cite{Matus62} (see also Bingham, Goldies and Teugel \cite{BiGoTe} Chapter 1 Theorem 1.5.4 p 23) that asserts the following: A nonnegative measurable function $L$ is slowly varying at $\infty$ iff for every 
$c\in (0, \infty)$ there exists a non-decreasing function $f_1$ and a non-increasing function $f_2$ such that  
$u^c L (u) \sim_{ \infty} f_1 (u)$ and $u^{-c} L (u) \sim_{ \infty} f_2 (u) $.  \cqfd

\medskip

To complete this section, we show that $\delta_\psi  >1$ is a more restrictive assumption than $\gamma_\psi >1$ by providing examples of branching mechanisms $\psi$ of the form (\ref{LevyKhin}), that satisfy (\ref{extinction}) and $1= \delta_\psi < \gamma_\psi$. 
\begin{lemma}
\label{contreex} For any $\gamma \in (1, 2]$, there exists a branching mechanism $\psi $ of the form (\ref{LevyKhin}) and such that $\eta_\psi= \gamma_{\psi}= \gamma$ and  $\delta_\psi= 1$. 
\end{lemma}
\noi
{\bf Proof:} For any $n \geq 3$, we set $\theta_n = n \log n$ and $\Delta_n = \theta_{n +1}-\theta_n$. It is easy to prove that $\Delta_n = \log n + 1 + {\cal O}( \frac{_1}{^n} ) \sim_{\infty}   \log n$. {\it We first suppose that $\gamma \in (1, 2)$}. For any $n \geq 3$, we set 
$r_n = \exp (- \theta_n)$, $a_n = r_n^{-\gamma}$ and $\pi(dr)= \sum_{^{n \geq 3}} a_n \delta_{r_n}(dr)$. It is easy to check that $\pi ((1, \infty))= 0$ and that $ \int_{(0, 1)} r^2 \pi (dr) = \sum_{^{n \geq 3}} 
r_n^{2-\gamma}$ is finite. We next define $\psi $ by 
\begin{equation}
\label{psipsi}
 \psi (\lambda ) = \int_{(0, \infty)} (e^{-\lambda r }-1+ \lambda r)    \, \pi (dr) \; , \quad \lambda \geq 0, 
 \end{equation}
that is clearly a branching mechanism of the form (\ref{LevyKhin}). We first prove that $\eta_\psi =\gamma_{\psi} =\gamma$, which is equivalent to $\eta_{\psi^\prime} =\gamma_{\psi^\prime} =\gamma -1$. Let us first prove that $\gamma_{\psi^\prime} \geq \gamma-1$. Note that 
$ \psi^\prime (\lambda ) = \sum_{^{n \geq 3} }   r_n^{-(\gamma-1)} (1-e^{-\lambda r_n}) $ 
and observe that for any $\lambda \geq 1/ r_3$, there exists $n_\lambda \geq 3$ such that $ r_{n_\lambda}^{-1} \leq \lambda < r_{n_\lambda + 1}^{-1}$. This inequality easily entails $n_\lambda < \log \lambda$. For all sufficiently large $\lambda $, we then get 
\begin{eqnarray*}
 \psi^\prime (\lambda) \geq r_{n_\lambda}^{-(\gamma-1)} (1-e^{-\lambda r_{n_\lambda}}) &\geq & (1-e^{-1}) e^{-(\gamma-1)\Delta_{n_\lambda}} \lambda^{\gamma-1} \\
 & \geq & (1-e^{-1}) \lambda^{\gamma-1} ( \log  \lambda )^{-2(\gamma-1)}  \; , 
 \end{eqnarray*} 
 which entails that $\gamma_{\psi^\prime} \geq \gamma-1$.

 Let us prove now that $\eta_{\psi^\prime} \leq \gamma-1$. To that end, we introduce the following notation 
 $$ R_n = \sum_{m \geq n} r_{m}^{2-\gamma} \quad {\rm and} \quad   S_n = \sum_{3 \leq m \leq n} r_{m}^{-(\gamma-1)} . $$
Elementary estimates entail that there exist two sequences $(\varepsilon_n , n \geq 0)$ and  
$(\varepsilon^\prime_n , n \geq 0)$, both converging to $0$, such that 
\begin{equation}
\label{estimRSn}
R_n = r_n^{2-\gamma} (1+ \varepsilon_n)  \quad{\rm and} \quad S_n = r_n^{-(\gamma-1)} (1+ \varepsilon^\prime_n) . 
\end{equation}
Since $ r_{n_\lambda}^{-1} \leq \lambda < r_{n_\lambda + 1}^{-1}$, we get 
\begin{eqnarray*}
 \psi^\prime (\lambda) & \leq & \!\!\!  \sum_{^{3 \leq n \leq n_\lambda} } \!\!\!  r_{n}^{-(\gamma-1)} +  \!\!\!\sum_{^{n \geq n_\lambda +1}}  \!\!\!  r_{n}^{-(\gamma-1)}  \lambda r_{n} \; \leq  S_{n_\lambda } +  \lambda R_{n_\lambda +1} \\
 & \leq & r_{n_\lambda }^{-(\gamma-1)} (1+  \varepsilon^\prime_{n_\lambda+1} ) + 
 \lambda r_{n_\lambda+1}^{2-\gamma}  (1+  \varepsilon_{n_\lambda +1} ) \\
 & \leq & \lambda^{\gamma -1} (2 + \varepsilon^\prime_{n_\lambda} +\varepsilon_{n_\lambda+1 }  ) \; , 
\end{eqnarray*} 
which shows that $\eta_{\psi^\prime} \leq \gamma-1$. We thus have proved $\eta_\psi=\gamma_\psi = \gamma$.

\smallskip

  Let us next prove that $\delta_\psi = 1$. We argue by contradiction and we suppose that $\delta_\psi >1$. Then, $\delta_{\psi^\prime}= \delta_\psi -1 >0$, and by (\ref{Jdeltaexpon}), there exist $c \in (0, 
  \delta_{\psi^\prime})$ and $C \in (0, \infty)$ such that 
\begin{equation}  
\label{explicitdelta}
C \,  y^{c-1} J_{\psi^\prime} (y) \leq  x^{c-1} J_{\psi^\prime} (x) \; , \quad 0 < x\leq y \leq 1 \; .
\end{equation}
Recall (\ref{Jdef}) and observe that $ J_{\psi^\prime } (x)= R_{n+1}+ x S_n$, for any $x\in [r_{n+1} , r_n) $. We set 
$$x_n = \exp \big( \!-\!(2-\gamma) \, \theta_{n+1} - \!(\gamma-1)\, \theta_{n} \big) = r_{n+1}^{2-\gamma} r_n^{\gamma-1} \; \in \,  [ \,  r_{n+1} \, , r_n \, )\;  . $$ Observe that 
$$ J_{\psi^\prime} (x_n)=   r_{n+1}^{2-\gamma} (2+ \varepsilon_{n+1} + \varepsilon^\prime_n )  \sim_{\infty} 2e^{-(2-\gamma)\theta_{n+1}}  \; .$$
Note that $ J_{\psi^\prime} (r_n) = R_n + r_n S_{n-1}$. This entails 
\begin{eqnarray*}
 J_{\psi^\prime} (r_n) &=&  e^{-(2-\gamma)\theta_n}  (1+ \varepsilon_n ) + 
e^{-\theta_n} e^{(\gamma-1)\theta_{n-1}}( 1+ \varepsilon^\prime_{n-1}) \\
& =& e^{-(2-\gamma)\theta_n} \big( 1 + \varepsilon_n + e^{-(\gamma-1)\Delta_{n-1}}( 1+ \varepsilon^\prime_{n-1}) \big) \\
& \sim_{\infty}  & e^{-(2-\gamma)\theta_n} \; .
\end{eqnarray*}
Thus, $ r_n^{c-1}J_{\psi^\prime} (r_n) \sim_\infty e^{ (\gamma-1-c)\theta_n}$ and 
$$  x_n^{c-1}J_{\psi^\prime} (x_n) \sim_\infty 2 e^{(1-c)(\gamma-1) \theta_n -c(2-\gamma) \theta_{n+1}}= e^{-c(2-\gamma) \Delta_n} e^{ (\gamma-1-c)\theta_n}  \; . $$
Recall that $\Delta_n \sim_\infty \log n$.  Since $0 <x_n < r_n \leq 1$, (\ref{explicitdelta}) imply that 
$$ 0 < \,  C \, \leq  \; \frac{x_n^{c-1}J_{\psi^\prime} (x_n)}{r_n^{c-1}J_{\psi^\prime} (r_n)} \; \sim_ \infty  \, e^{-c(2-\gamma) \Delta_n} \underset{n \rightarrow \infty}{ \longrightarrow} 0 , $$
which is aburd. This completes the proof of the lemma when $\gamma \in (1, 2)$. 

Let us consider the case $\gamma = 2$. For any $n \geq 2$, we set $r_n = e^{-n^2}$ and $\pi (dr)= \sum_{^{n \geq 2}} r_n^{-2}e^{-n\log n} \delta_{r_n} (dr)$. We define $\psi$ by (\ref{psipsi}) that is easily shown to be a branching mechanism of the form (\ref{LevyKhin}). We want to prove that $\gamma_{\psi}= 2$, (which implies that $\eta_\psi = 2$) and that $\delta_{\psi}=1$. Recall that it is equivalent to prove that 
$\gamma_{\psi^\prime}= 1$ and $\delta_{\psi^\prime}= 0$. For any $\lambda >e^4$, there exists
an integer $n_\lambda \geq 2 $ such that $e^{n_\lambda^2} \leq \lambda < e^{(n_\lambda +1)^2}$. Namely, $n_\lambda$ is the integer part of $ \sqrt{ \log \lambda }$. Observe that $\lambda r_{n_\lambda +1} <1$ and use the inequality $1-e^{-x} \geq x/2 $, $x \in [0, 1]$, to get the following inequality that holds true for all sufficiently large $\lambda$: 
\begin{eqnarray*}
 \psi^\prime (\lambda)& \geq & r_{n_\lambda +1}^{-1} e^{-(n_\lambda +1) \log (n_\lambda +1) } (1-
e^{ -\lambda r_{n_\lambda +1}} )   \geq  \frac{_1}{^2} \lambda e^{-(n_\lambda +1) \log (n_\lambda +1)} \\ 
&\geq & \lambda \exp (-2\sqrt{\log \lambda} \log \! \log \lambda ) , 
\end{eqnarray*}
which easily entails $\gamma_{\psi^\prime} \geq 1$, and thus $\gamma_{\psi^\prime} =\eta_{\psi^\prime} =1$, since $\psi^\prime$ is concave. We next set 
$$ R_n = \sum_{m \geq n} e^{-m\log m}  \sim_\infty e^{-n \log n}\quad {\rm and} \quad   S_n = \sum_{2 \leq m \leq n} e^{m^2 -m\log m} \sim_{\infty} e^{n^2 -n\log n}  . $$
Recall (\ref{Jdef}) and observe that $ J_{\psi^\prime } (x)= R_{n+1}+ x S_n$, for any $x\in [r_{n+1} , r_n) $. Recall notation $\Delta_n = (n+1)\log (n+1)-n\log n \sim_\infty \log n $. We next set $x_n: = r_ne^{-\Delta_n}$ that belongs to the interval $(r_{n+1}, r_n)$ for all sufficiently large integers $n$. It easy to check that for any $c \in (0, 1)$, one has   
$$ x_n^{c-1}J_{\psi^\prime} (x_n) \sim_\infty 2 e^{(1-c)n^2 -n \log n  -c \Delta_{n} }  \quad {\rm and} \quad r_n^{c-1} J_{\psi^\prime} (r_n) \sim_\infty e^{(1-c) n^2 -n\log n  }  $$
Thus, for any $c \in (0, 1)$, $\lim_\infty x_n^{c-1}J_{\psi^\prime} (x_n)/ r_n^{c-1} J_{\psi^\prime} (r_n) = 0$. This proves $\delta_{\psi^\prime} = 0$, which completes the proof of the lemma. \cqfd 

\subsection{Estimates.}
\label{estimsec}
In this section we state the estimates used in the proofs of Theorems \ref{packingLevytreeth} and \ref{densityth}. {\it Throughout the section we assume that $\psi$ is a branching mechanism of the form (\ref{LevyKhin}) whose exponent $ \delta$ defined by (\ref{doublingexponent}) is strictly larger than $1$}. Recall that $\varphi= \psi^\prime \circ \psi^{-1}$ and that 
$\varphi^{-1}$ stands for its reciprocal. Recall from (\ref{vvvequa}) the definition of the decreasing function  $v : (0, \infty) \rightarrow (0, \infty)$. 
\begin{lemma}
\label{controlvg} There exist $r_1 , C_1 \in (0, \infty)$, that only depend on $\psi$ and that satisfy 
$$ \forall r \in (0, r_1) \; , \quad   v(r) < C_1  r \varphi^{-1} (1/r) \; .$$ 
\end{lemma}
\noi
{\bf Proof:} Since $\delta >1$, there exist $c \in (1, \infty)$ and $C \in (0, \infty)$ such that $\psi (\lambda) \leq  C \psi(\lambda u)u^{-c}$, for any $u , \lambda \in [1, \infty)$. Choose $Q$ such that $C. \int^\infty_{Q} u^{-c} du \leq 1/4$. Thus,  
$$ \frac{\psi (\lambda)}{\lambda} \int_{\lambda Q}^\infty \frac{du}{\psi (u)}= \psi (\lambda) \int_{Q}^\infty \frac{du}{\psi (\lambda u)} \leq C. \int^\infty_{Q} u^{-c} du \leq \frac{1}{4} . $$
Denote by $v^{-1}$ the reciprocal of $v$ and recall that $v$ satisfies  (\ref{vvvequa}). Then, the previous inequality entails that $v^{-1} (Q\lambda) \leq \lambda / (4 \psi (\lambda) )$ and (\ref{convexpsi}) implies 
$v^{-1} (Q\lambda) \leq 1/\psi^\prime (\lambda)$. Since $v$ is decreasing we get $ v ( 1/ \psi^\prime (\lambda) ) \leq Q \lambda$. Substitute $\lambda $ with $\psi^{-1} (\lambda)$ to get 
$$ v \big( 1/ \varphi (\lambda) \big) \leq Q \psi^{-1} (\lambda ) \; , \quad \lambda \geq \psi (1) . $$
Next observe that $\psi^{-1} (\lambda ) \leq 4\lambda / \varphi (\lambda)$ by (\ref{convexpsi}). Thus, 
\begin{equation}
\label{controlvgconv}
v \big( 1/ \varphi (\lambda) \big) \leq \frac{4Q\lambda}{\varphi (\lambda)} \; , \quad \lambda \geq \psi (1) . 
\end{equation}
Set $C_1= 4Q$, $r_1= \varphi (\psi (1)) = \psi^\prime (1)$ and apply (\ref{controlvgconv}) with $\lambda = \varphi^{-1} (1/r)$ to get the desired result. \cqfd 

\smallskip

Recall from (\ref{kappadef}) the definition of $\kappa_r (\lambda, \mu)$ and recall that it satisfies the differential equation (\ref{equakappa}). Recall from (\ref{Mstardef}) the definition of $M_r^*(a)$. Observe that if $a \geq r$, then $M^*_r (a)= M^*_r (r)$. To simplify notation, we set 
\begin{equation}
\label{Lrlambdadef}
M^*_r := M^*_r  (r) \quad {\rm and} \quad  \cL_r (\lambda):= 1 - \frac{\psi \big( \kappa_r ( \lambda , 0)\big)}{\lambda} \;  , \quad r, \lambda \geq 0 \;. 
\end{equation}
\begin{lemma}
\label{Mstarrlaw} For any $ r \in (0, \infty)$, and for any $\lambda \in [0, \infty)$, one has 
 $$\cL_r (\lambda)  =N \left( \langle \ell^r\rangle e^{-\lambda \bm (\bar{B}(\rho, r))} \right) = e^{-\alpha r} \bE \left[ e^{-\lambda M^*_r }   \right]  \; .$$
\end{lemma}
\noi
{\bf Proof:} First observe that if $ \mu = \psi^{^{_{-1}}} (\lambda)$, then (\ref{kappaexc}) entails that $\kappa_r (\lambda , \mu)= \psi^{^{_{-1}}} (\lambda)$, for any $r \geq 0$. If  $\mu \neq \psi^{^{_{-1}}} (\lambda)$, then (\ref{kappaexc}) can be rewritten as the following integral equation
\begin{equation}
\label{integkappa}
\int^{\kappa_r (\lambda, \mu)}_\mu \!\!\!  \frac{du}{\lambda -\psi (u)} = r \; , \quad r , \lambda , \mu \geq 0 \quad {\rm and} \quad  \mu \neq \psi^{^{_{-1}}} (\lambda)\; .
\end{equation}
Note that $\mu \langle \ell^r\rangle + \lambda \bm (\bar{B}(\rho , r))= \mu L^r_\zeta + \lambda \int_0^\zeta ds \un_{\{ H_s \leq  r \}} $, for any $r, \lambda , \mu \geq 0$. Then, (\ref{kappaexc}) entails 
$N ( 1- \exp (-\mu \langle \ell^r \rangle -\lambda \bm (\bar{B}(\rho , r)) \, ) \, ) = \kappa_r ( \lambda , \mu ) $. 
We differentiate this identity with respect to $\mu$ to get 
$$ N \left( \langle \ell^r\rangle e^{-\mu \langle \ell^r\rangle -\lambda \bm (\bar{B}(\rho , r))} \right)= \frac{\partial \kappa_r}{\partial \mu} (\lambda , \mu) = \frac{\lambda -\psi \big( \kappa_r (\lambda , \mu)\big)}{\lambda -\psi (\mu)} ,  $$
which implies the first equality by taking $\mu= 0$. 

\medskip

It remains to prove that $ e^{\alpha r}\cL_r (\lambda)= \bE [ \exp (-\lambda M^*_r )  ] $. To that end, recall that $U$ is a (conservative) subordinator defined on $(\Omega, \cF, \bP)$ with Laplace exponent $\psi^{*\prime}= \psi^\prime-\alpha $. Then (\ref{Mstardef}) and (\ref{conditNstar}) imply 
$$ \bE \left[ \left. \exp (-\lambda M^*_r ) \right| \, U \, \right] = \exp \Big( -\int_{[0,r] }\!\!\!\!\!  dU_s \, \kappa_{r-s} (\lambda, 0)  \Big). $$
We therefore get 
\begin{eqnarray*}
 \bE \left[ \exp \Big( -\int_{[0, r] } \!\!\!\!\!   dU_s \kappa_{r-s} (\lambda, 0)  \Big) \right] &=&  \exp \Big(\! -\!\! \int_{0}^r \!\! ds \, 
\psi^{*\prime} ( \kappa_{r-s} (\lambda, 0) ) \,  \Big)  \\
& =& \exp \Big( \alpha r \!- \!\! \int_0^r \!\! ds\,  
\psi^{\prime} ( \kappa_{s} (\lambda, 0) ) \Big)  . 
\end{eqnarray*}
Now recall that $\frac{\partial}{\partial s} \kappa_s (\lambda, 0)= \lambda -\psi (\kappa_s (\lambda , 0))$ and a simple change of variable gives
$$ \int_0^r ds \, 
\psi^{\prime}\big( \kappa_{s} (\lambda, 0) \big)= \log \lambda -\log \big( \lambda - \psi (\kappa_r (\lambda , 0)) \,  \big)\; , $$ 
which easily completes the proof of the lemma. \cqfd 
\begin{lemma}
\label{Lrlambdainteg} For any $r, \lambda \geq 0$, one has 
$$ \int_0^{ -\log  \cL_r (\lambda) }\!\!\!\!\!\!\!\!\!\!\!\!\! \frac{dx}{\varphi \big( \lambda (1\!- \!e^{-x}) \big)} = r . $$
\end{lemma}
\noi
{\bf Proof:} Recall (\ref{integkappa}) that asserts that $r= \int_0^{\kappa_r (\lambda , 0)} du (\lambda -\psi (u))^{-1}$. Set $v= \psi (u)$ and recall that the derivative of $\psi^{-1}$ is $1/ \varphi$. So, easy changes of variable entail
$$ r=\!\! \int_0^{\kappa_r (\lambda , 0)} \!\!\!\!\!\!\!\!\!    \frac{du }{\lambda \!-\!\psi (u)} = \!\!\int_0^{\psi (\kappa_r (\lambda , 0) )}\!\!\!\!\!\!\!\!\!\!\! \frac{dv}{(\lambda \!- \!v)\varphi(v) } = \!\!\int_{\cL_r (\lambda)}^1
 \frac{w^{-1}dw}{\varphi( \lambda (1\!- \!w))}=\!\!\int_0^{ -\log  \cL_r (\lambda) }\!\!  \!\!\!\!\!\!\!\!\!\!\!\!\!\frac{dx}{\varphi \big( \lambda (1\!- \!e^{-x}) \big)}  $$
that is the desired result. \cqfd 
\begin{remark}
\label{Lrlambdainc} It is obvious from Lemma \ref{Mstarrlaw} that $\lambda \mapsto -\log \cL_r (\lambda)$ is increasing. Note that Lemma \ref{Lrlambdainteg} implies that $r \mapsto -\log \cL_r (\lambda)$ is also increasing. \cq 
\end{remark}
We now prove the key estimate for the lower bound in Theorem \ref{densityth}. 
\begin{lemma}
\label{densthkeyestim}  
Set $C_2 = (1-e^{-1})^{-1}$. There exists $r _2 \in (0, \infty)$ that only depends on $\psi$ such that 
$$ \cL_{2r} \big( C_2  \varphi^{-1} \big( \frac{_2}{^r}\log\! \log \frac{_2}{^r} \big) \,   \big) \; \leq  \; \exp \big(\!\!-2 \!\log \!\log\frac{_2}{^r} \, \big) \; , \quad  r \in (0, r_2). $$
\end{lemma}
\noi
{\bf Proof:} The proof is in four steps. We first claim the following. 
\begin{equation} 
\tag{{\rm Claim 1}}
\forall \; r, \lambda \in (0, \infty) , \quad  -\log \cL_{r} (\lambda ) \leq 1 \quad  \Longrightarrow \quad   \frac{2}{\lambda}\psi \big(  r \lambda / 2 \big)    \leq 1  \; \, . 
\end {equation}
\noi
{\it Proof of (Claim 1):} Note that $1-e^{-x} \geq x/2$ for any $x\in [0, 1]$ and recall that $1/\varphi$ is the derivative of $\psi^{-1}$. If $ -\log \cL_{r} (\lambda ) \leq 1$, then Lemma \ref{Lrlambdainteg} entails that 
$$ r = \int_0^{ -\log  \cL_r (\lambda) }\!\!  \!\!\!\!\!\!\!\!\!\!\!\!\!\frac{dx}{\varphi \big( \lambda (1\!- \!e^{-x}) \big)}
\leq \int_0^1 \frac{dx}{\varphi(\lambda x/2) }= 2 \psi^{-1} (\lambda /2) / \lambda \; , $$
which entails (Claim 1). We next claim the following. 
\begin{equation} 
\tag{{\rm Claim 2}}
\forall \; r, \lambda \in (0, \infty) , \quad  -\log \cL_{r} (\lambda ) >1   \Longrightarrow   \cL_{2r} (2 \lambda)  \leq  \exp \big( \!\!- \! r \, \varphi \big( \, 2(1\!-\!e^{_{-1}})  \lambda \, \big) \,  \big)   . 
\end {equation}
\noi
{\it Proof of (Claim 2):} Assume that $ -\log \cL_{r} (\lambda ) >1$. Then, Lemma \ref{Lrlambdainteg},  combined with elementary inequalities entails the following. 
\begin{equation}
\label{etapeclaim2}
\frac{\psi^{-1} (\lambda)}{\lambda} = \int_0^1 \frac{dx}{\varphi (\lambda x)} \leq \int_0^1 \frac{dx}{\varphi (\lambda(1\!- \!e^{- x}))} \leq \int_0^{ -\log  \cL_r (\lambda) }\!\!  \!\!\!\!\!\!\!\!\!\!\!\!\!\frac{dx}{\varphi \big( \lambda (1\!- \!e^{-x}) \big)}= r . 
\end{equation}
By Remark \ref{Lrlambdainc}, we have $-\log \cL_{2r} (2\lambda ) >-\log \cL_{r} (\lambda ) >1$. Thus, we get 
\begin{eqnarray}
\label{swippage}
2 r &=& \int_0^{ -\log  \cL_{2r} (2\lambda) }\!\!  \!\!\!\!\!\!\!\!\!\!\!\!\!\frac{dx}{\varphi \big( 2\lambda (1\!- \!e^{-x}) \big)} \leq \int_0^1 \frac{dx}{\varphi (\lambda x)} + \int_1^{ -\log  \cL_{2r} (2\lambda) }\!\!  \!\!\!\!\!\!\!\!\!\!\!\!\!\frac{dx}{\varphi \big( 2\lambda (1\!- \!e^{-x}) \big)} \nonumber \\ 
& \leq & \frac{\psi^{-1} (\lambda)}{\lambda}  - \frac{\log  \cL_{2r} (2\lambda)}{\varphi \big( 2 (1\!- \!e^{-1}) \lambda \big)} . 
\end{eqnarray} 
(here again, we use the inequality $1-e^{-x} \geq x/2$, $x \in [0, 1]$). Then (\ref{swippage}) and (\ref{etapeclaim2}) entail that
$r \leq -\log  \cL_{2r} (2\lambda)/ \varphi ( 2 (1\!- \!e^{-1}) \lambda )$, which implies (Claim 2).

\smallskip

Recall from (\ref{gaugedef}) the definition of $g$. We claim that there exists $R \in (0, 2 r_0)$ such that 
\begin{equation}
\tag{{\rm Claim 3}}
\forall r \in (0, R) \; , \quad -\log \cL_{r} \Big( \frac{_4}{^{g(r/2)}} \Big) >1 \; .
\end{equation}
\noi
{\it Proof of (Claim 3)}: Let us set $ \lambda_{r} = \varphi^{-1} ( \frac{2}{r}\log\! \log \frac{2}{r})$, for any $r \in (0, r_0)$. Thus, $g(r/2)= (\log\! \log \frac{2}{r})/ \lambda_r$. Suppose that $-\log \cL_{r} (\frac{4}{ g(r/2) }) \leq 1$. Then, (Claim 1) easily entails that $ 1 \geq   \frac{1}{2} g(r/2) \psi \big( 4r/(2g(r/2)) \, \big)$, which is equivalent to the following: 
$$ 1 \geq \frac{1}{2 }\cdot \frac{\log\! \log \frac{2}{r} }{ \lambda_{r}} \cdot \psi 
\left(  \frac{4\lambda_{r}}{\frac{2}{r} \log\! \log \frac{2}{r} } \right)= \frac{1}{2 }\cdot \frac{\log\! \log \frac{2}{r} }{ \lambda_{r}} \cdot \psi 
\left(  \frac{ 4\lambda_{r}}{\varphi(\lambda_{r}) } \right) $$
Now recall from (\ref{convexpsi}) that $4 \lambda_{r}/ \varphi(\lambda_{r}) \geq  \psi^{-1} (\lambda_{r}) $. Thus, the latter inequality implies $ 2 \geq  \log\! \log \frac{2}{r} $ and (Claim 3) holds true with $R$ being the largest $r $ in $(0, 2r_0) $ such that $\log\! \log \frac{2}{r}  \geq  2$. 
\medskip

\noi
{\it End of the proof of the lemma:} Recall the notation $C_2=(1-e^{-1})^{-1}$. There exists $r_2 \in (0, R)$ such that $ 4/ \log\! \log \frac{2}{r_2} \leq C_2/2$. The definition of $g$ implies 
$$ \frac{4}{g(r/2)}= \frac{4}{\log\! \log \frac{2}{r}}  \, \varphi^{-1} \big( \frac{_2}{^r}\log\! \log \frac{_2}{^r} \big) \leq \frac{_1}{^2} C_2\varphi^{-1} \big( \frac{_2}{^r}\log\! \log \frac{_2}{^r} \big) \; , \quad r \in (0, r_2) . $$
Remark \ref{Lrlambdainc} and (Claim 3) entail that for any $r \in (0, r_2)$, 
$$ 1< -\log \cL_{r} \Big( \frac{_4}{^{g(r/2)}} \Big) \leq - \log \cL_r \Big( \frac{_1}{^2} C_2 \varphi^{-1} \big( \frac{_2}{^r}\log\! \log \frac{_2}{^r} \big) \Big) . $$
Now (Claim 2) implies that 
$$ \cL_{2r} \left( C_2 \varphi^{-1} \big( \frac{_2}{^r}\log\! \log \frac{_2}{^r} \big) \right) \leq  \exp \left( -r \varphi \big( \varphi^{-1} \big( \frac{_2}{^r}\log\! \log \frac{_2}{^r} \big) \big) \right) =
 \exp \left( -2 \log\! \log \frac{_2}{^r} \right), $$
which completes the proof of the lemma. \cqfd 

\smallskip

The following estimate is used in the proof of Theorem \ref{packingLevytreeth}. 
\begin{lemma} 
\label{treillisestim} There exist $r_3, C_3, C_4 \in (0, \infty)$ that only depend on $\psi$ such that for any $r \in (0, r_3)$, one has 
$$ g(16r) \, N \Big( \sup H \geq 3r  \, ;  \!\int_0^\zeta \!\!\!\! \un_{\{ H_s \leq 2r  \}} ds  \leq  C_3\,  g(16  r)  \Big) \leq C_4 e^{-  \frac{3}{2} \log \! \log \frac{2}{r}}  \,  r  \log \!\log \frac{_1}{^r}   .$$\end{lemma}
\noi
{\bf Proof:} Recall that $(g_j^{2r}, d^{2r}_j) $, $j \in \cI_{2r}$, stand for the open connected components of 
$\{ s \in [0, \zeta ] \; : H_{s} > 2r \}$ and that $H^{2r,j}= H_{(g_j^{2r} + \cdot ) \wedge d^{2r}_j}-2r$, $j \in \cI_{2r}$, are the corresponding excursions of $H$ above $2r$. Recall that $\cG_{2r}$ is the sigma-field generated by the height process $\widetilde{H}^{2r}$ below level $2r$, augmented by the $N$-negligible sets. Recall from (\ref{condiprob}) the notation $N_{2r}$.  
The branching property asserts that under $N_{2r}$ and conditionally given  $\cG_{2r}$, the random variable $ Y:=  \sum_{^{j\in \cI_{2r}}} \un_{ \{  \sup H^{2r,j}  \geq r  \}}$
is distributed as a Poisson random variable with parameter $L^{2r}_\zeta N ( \sup H >r)= L^{2r}_\zeta v(r)$. Now observe that $N$-a.e. $\un_{\{ Y \neq 0\}}= \un_{\{ \sup H \geq 3r\}}$. Thus, 
$$ N_{2r} \left( \un_{\{ \sup H \geq 3r\}} \, | \, \cG_{2r} \right)= N_{2r} \left( \un_{\{ Y \neq 0\}} \, | \, \cG_{2r} \right)= 1-e^{-v(r)L^{2r}_\zeta} \leq v(r) L^{2r}_\zeta . $$
Since $ \int_0^\zeta \un_{\{ H_s \leq 2r  \}} ds$ is $\cG_{2r}$-measurable, we get 
$$N_{2r} \Big( \sup H \geq 3r  \, ;  \!\int_0^\zeta \!\!\!\! \un_{\{ H_s \leq 2r  \}} ds  \leq  C_3 g(16r)  \Big) \leq 
v(r) N_{2r}\big( L^{2r}_\zeta \un_{ \{ \!\int_0^\zeta \!  \un_{\{ H_s \leq 2r  \}} ds  \leq  C_3  g(16 r) \}  } \big). $$
where $C_3$ is a positive constant to be specified further. Recall that $N$-a.e. $\un_{\{ L^{2r}_\zeta \neq 0 \}}= \un_{\{ \sup H >2r \}}$. Consequently, 
$$ N \Big( \sup H \geq 3r  \, ;  \!\int_0^\zeta \!\!\!\! \un_{\{ H_s \leq 2r  \}} ds  \leq  C_3 g(16r)  \Big) \leq 
v(r) N\big( L^{2r}_\zeta \un_{ \{ \!\int_0^\zeta \!  \un_{\{ H_s \leq 2r  \}} ds  \leq  C_3 g(16r) \}  } \big). $$
Recall that $L^{2r}_\zeta= \langle \ell^{2r} \rangle$, that $\int_0^\zeta \!  \un_{\{ H_s \leq 2r  \}} ds= \bm (\bar{B} (\rho , 2r))$ and recall  from (\ref{Lrlambdadef}) the notation $\cL_{2r} (\lambda)$. Then, the Markov inequality combined with Lemma \ref{Mstarrlaw} entails for any $\lambda \geq 0$, 
\begin{eqnarray}
\label{logpart}
N\big( L^{2r}_\zeta \un_{ \{ \!\int_0^\zeta \!  \un_{\{ H_s \leq 2r  \}} ds  \leq  C_3 g(16r) \}  } \big)  &= &
 N \big( \langle \ell^{2r} \rangle \un_{\{ \bm (\bar{B} (\rho , 2r)) \leq  C_3 g(16r)  \}} \big)  \nonumber \\ 
& \leq & e^{C_3 \lambda g(16r)} \cL_{2r} \big( \lambda  \big) .
\end{eqnarray}
Set $r_3= r_1 \wedge r_2$, where  $r_1$ and $r_2$ are as in Lemmas \ref{controlvg} and \ref{densthkeyestim}. We fix $r \in (0, r_3)$. Since we assumed that $\delta =\delta_\psi >1$, there exists $C\geq 1$ such that $g$ satisfies a $C$-doubling condition (\ref{doubling}). Thus, $g( 16r) \leq C^5 g(r/2)$. Recall notation $\lambda_r = \varphi^{-1} ( \frac{2}{r} \log \! \log \frac{2}{r})$ and note that $\lambda_r g(r/2)= \log \! \log \frac{2}{r}$. Take $\lambda= C_2 \lambda_r$ in (\ref{logpart}) and use Lemma \ref{densthkeyestim} to get 
$$  N\big( L^{2r}_\zeta \un_{ \{ \!\int_0^\zeta \!  \un_{\{ H_s \leq 2r  \}} ds  \leq  C_3 g(16r) \}  } \big)  \leq 
\exp \big( -(2- C_3 C^5) \log \! \log \frac{_2}{^r} \big). $$
We now choose $C_3 $ such that $C_3C^5= 1/2$. We then get for any $r \in (0, r_3)$, 
\begin{equation}
\label{etapepepe}
g(16r) N \Big( \sup H \geq 3r  \, ;  \!\int_0^\zeta \!\!\!\! \un_{\{ H_s \leq 2r  \}} ds  \leq  C_3 g(16r)  \Big) \leq  g(16r) v(r)  e^{-  \frac{3}{2} \log \! \log \frac{2}{r}} . 
\end{equation}
We now use Lemma \ref{controlvg} to get the following.  
\begin{eqnarray*}
g(16r) v(r) & \leq & C^4 \cdot g(r) v(r) \leq  C^4 \cdot \frac{v(r) \log \! \log \frac{1}{r}}{\varphi^{-1} ( \frac{1}{r } \log \! \log \frac{1}{r} ) } \\
 & \leq & C^4 \cdot \frac{v(r) \log \! \log \frac{1}{r}}{\varphi^{-1} ( \frac{1}{r } ) }  \leq  C^4 C_1 \cdot r \log \! \log \frac{_1}{^r} \; , 
\end{eqnarray*}
which implies the desired result with $C_4= C^4C_1$ by (\ref{etapepepe}). \cqfd 

\medskip

We recall in a lemma a result due to Fristedt and Pruitt \cite{FriPru71} on subordinators that is needed to prove the upper bound in Theorem \ref{densityth}.  Recall that $\varphi= \psi^\prime \circ \psi^{-1}$ and that $\psi^\prime (0)= \alpha$. We set $\varphi^* = \varphi -\alpha$. Thus, $\varphi^* (0)= 0$ and $\varphi^*$ is the Laplace exponent of a conservative subordinator.  
\begin{lemma}
\label{subordpropri} (Fristedt and Pruitt \cite{FriPru71} Theorem 1 p 173) Let $(S_r, r\geq 0)$ be a subordinator starting at $0$ defined on $(\Omega, \cF, \bP)$ whose Laplace exponent is $\varphi^*$.  Let us assume that $\delta =\delta_{\psi} >1$. Recall that $g$ is the gauge function defined by (\ref{gaugedef}). Then, there exists a constant $K_{\psi} \in (0, \infty)$ that only depends on $\psi$ such that 
$$ \bP^{_{_-}}{\rm a.s.} \qquad  \liminf_{r\rightarrow 0} \frac{S_r}{g(r)} = K_\psi \; .$$
\end{lemma}

\begin{remark} 
\label{aegalun}  Theorem 1 \cite{FriPru71} is actually more general than the result stated in Lemma \ref{subordpropri}. It actually asserts the following:  If  
$\gamma= \gamma_\psi >1$, then for any $a >1$,  there exists $K_{a, \psi} \in (0, \infty)$ that only depends on $a$ and $\psi$,  such that 
\begin{equation}
\label{vraiFriPru}
 \bP^{_{_-}}{\rm a.s.} \qquad  \liminf_{r\rightarrow 0} \frac{S_r}{g^*(ar)} = K_{a, \psi} \; , 
 \end{equation}
where $g^*$ is derived from $\varphi^* $ as $g$ is derived from $\varphi$ in (\ref{gaugedef}). Observe that $\varphi^* $ and $\varphi$ are equivalent at infinity. Thus, $g \sim_0 g^*$ and $g^*$ can be replaced by $g$ in (\ref{vraiFriPru}). If $\delta =\delta_{\psi} >1$, then $g$ satisfies the doubling condition (\ref{doubling}). Consequently, $g(ar)\asymp_0 g(r)$ and Blumenthal zero-one law allows to deduce Lemma \ref{subordpropri} from (\ref{vraiFriPru}) under the more restrictive 
assumption $\delta =\delta_{\psi} >1$. \cq
\end{remark}

\section{Proof of Theorem \ref{densityth}.}
\label{proofdensthsec}
Recall  that $(U_t, t \geq 0)$ is a subordinator defined on $(\Omega, \cF, \bP)$ with Laplace exponent $\psi^{*\prime} (\lambda )= \psi^\prime (\lambda)-\alpha $, $\lambda \geq 0$. Recall that $\cN^* = \sum_{^{j \in \cI^*}} \delta_{ (  r^*_j , \, H^{*j} )  } $
is a random point measure on $[0, \infty) \times \cC^0$ defined on  $(\Omega, \cF, \bP)$ such that conditionally given $U$, $\cN^*$ is distributed as a Poisson point measure with intensity $dU_r \otimes N(dH)$. Also recall the notation
$$M^*_r= M^*_r (r) = \sum_{j \in \cI^*} \un_{[0, r]} (r^*_j) \cdot \int_{0}^{\zeta^*_j} \un_{\{ H^{*j}_s \leq r-r^*_j  \}} ds  \; , \quad r \geq 0 \; , $$
where $\zeta^*_j$ stands for the lifetime of $H^{*j}$, for any $j \in \cI^*$. We now set 
$$ \forall r \geq 0 \; , \quad S_r= \sum_{j \in \cI^*} \un_{[0, r]} (r^*_j) \zeta^*_j . $$
First observe that 
\begin{equation}
\label{domina}
 \forall r \geq 0 \; , \quad M^*_r \leq S_r \; .
\end{equation}
Next observe that $(S_r, r \geq 0)$ is cadlag and that $S_0= 0$. Let $ 0=r_0  <  r_1< \ldots < r_n $ and $\lambda_1, \ldots , \lambda_n \geq 0$. Recall from (\ref{lifetimeexc}) that $N(1-\exp (-\lambda \zeta)\,)= \psi^{-1} (\lambda)$. We then get 
\begin{eqnarray*} 
\bE \big[ \exp \Big( \!- \!\! \sum_{^{1\leq k \leq n}} \!\!  \lambda_k (S_{r_k}-S_{r_{k-1}}) \, \Big) \Big] &=& \bE
 \big[ \exp \Big( \! - \!\! \sum_{^{1\leq k \leq n}} \!\! 
\psi^{-1}(\lambda_k) (U_{r_k}-U_{r_{k-1}} ) \, \Big) \Big] \\
& =&  \exp \Big( \!- \!\! \sum_{^{1\leq k \leq n}} \!\!  (r_k-r_{k-1})\varphi^{*}(\lambda_k) \,  \Big). 
\end{eqnarray*}
This implies that $(S_r, r \geq 0)$ is a subordinator with Laplace exponent $\varphi^*$. Lemma \ref{subordpropri} applies: There exists a constant $K_\psi \in (0, \infty)$ that only depends on $\psi$ such that
\begin{equation}
\label{upperdensth}
\bP^{_{_-}}{\rm a.s.} \qquad \liminf_{r \rightarrow 0} \; \frac{M^*_r}{g(r)} \; \leq \; \liminf_{r \rightarrow 0} \; \frac{S_r}{g(r)} \; = K_\psi < \infty \; . 
\end{equation}
Let us prove a lower bound. Lemma \ref{Mstarrlaw} asserts that $\bE [ \exp( -\lambda M^*_r) ]= e^{\alpha r} \cL_r (\lambda) $. 
Recall notation $\lambda_r = \varphi^{-1} ( \frac{2}{r} \log \! \log \frac{2}{r})$ so that $\lambda_r g(r/2)= \log \! \log \frac{2}{r}$. Lemma \ref{densthkeyestim} asserts that there exist $r_2 , C_2 \in (0, \infty)$ such that $\cL_{2r} (C_2 \lambda_r) \leq \exp(-2\log \! \log \frac{2}{r})$ for any $r \in (0, r_2)$. Let $K$ be a positive real number to be specified later. By Markov inequality, for any $r \in (0, r_2)$, one has 
\begin{eqnarray*}
\bP \left( M^*_{2r} \leq K  g(r/2)  \, \right) &=& \bP \left( C_2  \lambda_r M^*_{2r} \leq C_2 K \log \! \log \frac{_2}{^r}  \right) \\
& \leq & \exp \big(  C_2K \log \! \log \frac{_2}{^r} \big) \bE  \left[ \exp \big(   -C_2  \lambda_r M^*_{2r}\big) \right] \\
& \leq & \exp \big(  C_2K  \log \! \log \frac{_2}{^r}  + 2 \alpha r \big) \cL_{2r} (C_2 \lambda_r) \\
& \leq & \exp \big(  -(2-C_2 K)   \log \! \log \frac{_2}{^r}  + 2 \alpha r \big) .
\end{eqnarray*}
We choose $K$ such that $C_2K= 1/2$. Borel-Cantelli entails  
$$ \textrm{$\bP$-a.s.} \qquad \liminf_{n \rightarrow \infty}  \frac{M^*_{2^{-n}}}{ g(2^{-n-2})} \geq K \; .$$
Recall that $g$ satisfies a $C$-doubling condition. Thus, $g(2r) \leq C^2 g(r/2)$, for all sufficiently small $r >0$. Since $r \mapsto M^*_r $ is non-decreasing, we get 
\begin{equation}
\label{lowerdensth}
 \textrm{$\bP$-a.s.}\qquad \liminf_{r \rightarrow 0} \; \frac{M^*_{r}}{g(r)} \geq  C^{-3} \liminf_{n \rightarrow \infty} \; \frac{M^*_{2^{-n}}}{g(2^{-n-2})} \geq C^{-3}K >0  \; .
 \end{equation}
By standard results on Poisson point processes, $\liminf_{r \rightarrow 0} \; M^*_{r} /g(r)$ 
is a random variable that is measurable with respect to the tail sigma-field of $U$ at $0+$ and Blumenthal's 0-1 law entails that there exists a constant denoted by $C_\psi$ such that 
$\bP$-a.s.$\, \liminf_{r \rightarrow 0} \; M^*_{r} /g(r)= C_\psi$. Moreover (\ref{upperdensth}) and (\ref{lowerdensth}) show that $C_\psi \in (0, \infty)$. Then,  (\ref{keymass}) implies 
$$ N \left( \int_{\cT} \!\! \bm (d\sigma) \,  \un_{ \{ \liminf_{r \rightarrow 0}  \bm (\overline{B}(\sigma, r)) / g(r) \neq C_\psi    \}} \right) = 0 \; , $$
which entails Theorem \ref{densityth}.  \cqfd

\section{Proof of Theorem \ref{packingLevytreeth}.}
\label{proofpackingth}
The proof is in several steps that are stated as lemmas. Recall that for any $a >0$, the intervals $(g_j^{a}, d^{a}_j) $, $j \in \cI_{a}$, are the open connected components of 
$\{ s \in [0, \zeta ] \; : H_{s} > a \}$. Recall notation $H^{a,j}= H_{(g_j^{a} + \cdot ) \wedge d^{a}_j}-a$, $j \in \cI_{a}$, for the excursions of $H$ above $a$. We shall denote by $\zeta_{a,j}= d^{a}_j-g_j^{a}$, the lifetime of $H^{a,j}$. Recall that $\cG_{a}$ stands for sigma-field generated by the height process $\tilde{H}^{a}$ below $a$, augmented by the $N$-negligible sets and recall from (\ref{condiprob}) the notation $N_{a}$. The branching property asserts that under $N_a$ and conditionally given $\cG_{a}$, the point measure 
$$\cM_a =  \sum_{j\in \cI_{a}} \delta_{\big(L^a_{g^a_j}\;  , \;  H^{a, j} \big)}$$
is distributed as a Poisson point measure with intensity  $\un_{[0, L^a_\zeta]} (x) dx \otimes N (dH)$. We apply this property as follows. Let $F: \cC^0 \rightarrow [0, \infty)$ be measurable. Set $Z^{a}_{F}=   \sum_{^{j\in \cI_{a}}} F( H^{a, j})$ and note that $Z^a_F= 0$ if  $ \sup H \leq a$. A basic result on Poisson point processes entails that $N_a ( Z^a_F | \cG_a )= L^a_\zeta N( F(H))$. Recall that $N$-a.e.$\, \un_{\{ L^a_\zeta \neq 0 \}}= \un_{\{ \sup H >a \}}$. Thus, $N (Z^a_F )= N(L^a_\zeta)N(F(H))$. By (\ref{meanloc}), $N(L^a_\zeta)= \exp (- \alpha a) \leq 1$, which entails  
\begin{equation}
\label{compenn}
N \Big( \sum_{^{j\in \cI_{a}}}   F( H^{a, j}) \Big) \leq N(F(H)) \; .
\end{equation}

For any $n\in \bN$, we denote by $\cD_{n}$ the set $\{ k2^{-n-3} \, ; \, k \in \bN\}$. For any 
$a \in \cD_n$ and for any $j \in \cI_a$, we define the event $E(a, j)$ as follows
$$ E(a, j) =   \Big\{  \; \sup H^{a,j} > 3.2^{-n-3}  \; ;  \! \int_{0}^{\zeta_{a,j }} \!\!\!\!\! \un_{\{ H^{a,j}_s \leq 2.2^{-n-3} \}} ds  \leq C_3 g( 2.2^{-n})  \; \Big\} , $$
where $C_3$ is the constant appearing in Lemma \ref{treillisestim}. For any positive integer $Q $, we then set
\begin{equation}
\label{UnQdef}
 U_{n, Q}=   \sum_{a \in \cD_{n}\cap [0, Q +1]} \; \sum_{ j \in \cI_a } g(2.2^{-n}) \, \un_{ E(a,j)  } \; \; . 
\end{equation}
We suppose that $2^{-n-3} < r_3$, where $r_3$  is the constant appearing in Lemma \ref{treillisestim}. We apply (\ref{compenn}) with $a \in \cD_{n}\cap [0, Q +1]$ and Lemma \ref{treillisestim} with $r_n= 2^{-n-3}$, and we get 
$$ N \Big(  \sum_{^{j \in \cI_a}} g(2.2^{-n}) \, \un_{ E(a,j)  }  \Big) \leq C_4 r_n \log \! \log \frac{_1}{^{r_n}} e^{- \frac{_3}{^2}\log \! \log \frac{_1}{^{r_n}}   }  \leq C_5 2^{-n-3} n^{-3/2}\log n ,  $$
Here $C_4$ is the constant appearing in Lemma \ref{treillisestim} and $C_5$ only depends on $\psi$. 
Since $\# (\cD_{n}\cap [0, Q +1] )= 2^{n+3}(Q+1) \leq 2Q .2^{n+3}$, we get 
\begin{equation}
\label{ligneestim}
N \left(   U_{n, Q}  \right) \leq  2C_5 Q \cdot    n^{-3/2}\log n \; , 
\end{equation}
which easily entails the following lemma. 
\begin{lemma}
\label{treillisas} There exists a Borel set $A_1\subset \cC^0$ such that $\cC^0 \backslash A_1$ is $N$-negligible and such that on $A_1$, 
$$  \forall Q \in \bN \; , \quad \lim_{n \rightarrow \infty} \; \sum_{m \geq n} U_{m, Q} \; = \, 0 \; . $$
\end{lemma}
Recall notation $\cP^*_g$ for the $g$-pre-packing measure on $\cT$. Lemma \ref{treillisas} is now used to prove the following.
\begin{lemma}
\label{dominmass} There exists $C_6 \in (0 , \infty) $ that only depends on $\psi$ such that on $A_1$,  one has $ \cP^*_g (K) \leq C_6 \bm (K) $, for any compact subset $K$ of $\cT$. 
\end{lemma}
\noi
{\bf Proof:} Recall the notation $\Gamma= \sup H= \sup \{ d( \rho , \sigma ) \; ; \; \sigma \in \cT \}$ that is the total height of $\cT$. Without loss of generality, we can assume that $\Gamma \in (0, \infty)$ on $A_1$. Let $0 < \varepsilon <  \min (1, \Gamma )$ and let $Q $ be a positive integer such that $Q > \Gamma $. Let $K$ be any compact subset of $\cT$ and let $(\bar{B}(\sigma_\ell , r_\ell) , \ell \geq 0)$ be any $\varepsilon$-packing of $K$. Namely, the closed balls $\bar{B}(\sigma_\ell , r_\ell) $'s are pairwise disjoint, $\sigma_\ell \in K$ and $r_\ell \leq \varepsilon$, for any $\ell \geq 0$. Observe that at most {\it one} ball may contain the root $\rho$. Without loss of generality, we assume that if a ball of the $\varepsilon$-packing of $K$ contains $\rho$, then it is $\bar{B}(\sigma_0 , r_0)$. Thus, $\rho \notin \bar{B}(\sigma_\ell , r_\ell)$, for any $\ell \geq 1$. 

   We fix $(\sigma , r ) \in \{ (\sigma_\ell , r_\ell ) \; ; \; \ell \geq 1 \}$. Since $r \leq \varepsilon < 1$, there exists an integer $n(r) \geq 1$ such that $2^{-n (r)} < r \leq 2.2^{-n(r)}$. Since $\rho \notin \bar{B} (\sigma , r)$, one has $d(\rho, \sigma) >r$ and there exists $t_0 \in (0, \zeta ) $ such that $p (t_0)= \sigma $ (recall that $p$ stands for the canonical projection from $[0, \zeta]$ onto $\cT$). Thus, 
$H_{t_0} =   d(\rho, \sigma) >r > 8.2^{-n(r)-3} $ and there exists a unique integer $k > 8$ such that $ k2^{-n(r)-3} \leq H_{t_0} < (k +1)2^{-n(r)-3}$. We then set  
\begin{displaymath}
\left \{    \begin{array}{l} 
 {\rm g} =\sup \{ s \in [0, t_0]  : \, H_s = (k-3)2^{-n(r)-3} \} \\
{\rm d} =\inf \{ s \in [t_0, \zeta]  : \, H_s = (k-3)2^{-n(r)-3} \}.
\end{array} \right.
\end{displaymath}
To simplify notation, we set $a= (k-3)2^{-n(r)-3} \in \cD_{n(r) } \cap [0, Q+1] $. 
\begin{description}
\item{{\bf (I)}} Observe that $0 < {\rm g} < t_0< {\rm d}  < \zeta$, that $H_{{\rm g}} = H_{{\rm d}}= a$, that $H_s >a$ for any $s \in ({\rm g}, {\rm d})$ and that $H_{t_0} -a \geq  3.2^{-n(r)-3}$. Therefore there exists a {\it unique} $j \in \cI_a$ such that 
$$  (g^a_j, d^a_j)= ({\rm g}, {\rm d}) \, , \;  H^{a, j}_s= H_{({\rm g} +s)\wedge {\rm d}} -a\;  , \; s  \geq 0\, , \quad {\rm and} \quad \sup H^{a, j} \geq 3 . 2^{-n(r)-3} \; .$$
\item{{\bf (II)}} Recall that $b(s_1, s_2)= \inf_{[s_1 \wedge s_2 , s_1 \vee s_2 ]} H$ and $d(s_1, s_2) = H_{s_1} + H_{s_2} -2 b(s_1, s_2) $, for any $s_1, s_2 \in [0, \zeta ]$.  Let $s \in ({\rm g}, {\rm d})$ be such that $H_{s}-a \leq 2.2^{-n(r)-3}$. First observe that  $b(s,t_0) \geq a$. Next observe that $H_{t_0}-a \leq 4.2^{-n(r)-3}$. Consequently, $d(s,t_0) \leq 2.2^{-n(r)-3} + H_{t_0} -a \leq 6. 2^{-n(r)-3} < r$. Thus,  
$$ \{ s \in ({\rm g}, {\rm d}) \; : \; H_{s}-a \leq 2.2^{-n(r)-3} \}  \subset \{ s \in [0, \zeta] \; : \; d(s,t_0) \leq r \} , $$
which implies that $ \int_0^{\zeta_{a, j}}  \un_{ \{ H^{a, j}_s  \leq 2. 2^{-n(r)-3} \}}ds  \leq \bm \big( \bar{B} (\sigma , r) \big) $. 
\end{description}
Recall that $\bar{B}(\sigma_0 , r_0)$ is the only ball of the $\varepsilon$-packing that may contain the root. Since the $\bar{B}(\sigma_\ell , r_\ell)$'s are pairwise disjoint, ${\bf (I)}$ and ${\bf (II)}$ imply the following inequalities
$$
\sum_{\ell \geq 0} g(r_\ell) \un_{ \{  \bm ( \overline{B} (\sigma_\ell , r_\ell ) ) \leq C_3 g(r_\ell)  \}}  
\;  \leq  \; g( \varepsilon )  + \sum_{ n : 2^{-n} < \varepsilon } U_{n, Q} \; .
$$
This entails the following
\begin{eqnarray}
\label{control}
\sum_{\ell \geq 0} g(r_\ell)  &=& \sum_{\ell \geq 0} g(r_\ell) \un_{ \{  \bm ( \overline{B} (\sigma_\ell , r_\ell) )  >C_3 g(r_\ell)  \}} + \sum_{\ell \geq 0} g(r_\ell) \un_{ \{  \bm ( \overline{B} (\sigma_\ell , r_\ell) ) \leq C_3 g(r_\ell)  \}}  \nonumber \\
 \! \!\! &\leq &  \! \!\!  
 \frac{_1}{^{C_3}}\bm \big( K^{(\varepsilon)}\big) + g( \varepsilon ) \! + \!\! \sum_{ n : 2^{-n} < \varepsilon } U_{n, Q} \; , 
\end{eqnarray}
where we have set $ K^{(\varepsilon)}= \{ \sigma \in \cT : \, d( \sigma , K) \leq \varepsilon \} $. Recall from (\ref{prepremeadef}) the definition of $\cP^{(\varepsilon)}_g$. Since the previous inequality holds true for any $\varepsilon$-closed packing of $K$, we get 
$$ \cP^{(\varepsilon)}_g ( K) \leq  \frac{_1}{^{C_3}}  \bm \big( K^{(\varepsilon)}\big) + g( \varepsilon ) \! + \!\! \sum_{ n : 2^{-n} < \varepsilon } U_{n, Q} \; . $$
Since $K$ is compact, $\lim_{\varepsilon \rightarrow 0} \bm (K^{(\varepsilon)})= \bm (K)$. The previous inequality, combined with Lemma \ref{treillisas}, implies the desired result on $A_1$ with $C_6=1/C_3 $. \cqfd

\smallskip 

Recall that Theorem \ref{densityth} asserts that there exists a Borel subset $A_2 \subset \cC^0$ such that $N ( \cC^0 \backslash A_2)= 0$ and such that, on $A_2$, we have 
\begin{equation}
\label{badmass}
\bm \left(  \big\{ \sigma \in \cT \; : \; \liminf_{r \rightarrow 0} g(r)^{-1} \bm \big( B( \sigma , r) \big) \neq  
C_\psi \;   \big\} \right) = 0 \;  
\end{equation}

\begin{lemma} 
\label{masspackequiv}
There exist $C_7, C_{8} \in (0, \infty)$ that only depend on $\psi$ such that on $A_1\cap A_2$, for any Borel set $B\subset \cT$, we have 
$$ C_7 \bm (B) \leq \cP_g (B) \leq C_{8}  \bm (B) \; .$$
\end{lemma}
\noi
{\bf Proof:}  We set $ {\rm Good} =\{ \sigma \in \cT :  \liminf_{r \rightarrow 0} g(r)^{-1} \bm ( B( \sigma , r) ) = 
C_\psi  \}$ and ${\rm Bad}= \cT \backslash {\rm Good}$. We argue deterministically on $A_1\cap A_2$.  
Observe that $\cP_g ( {\rm Bad}) \leq \cP^*_g (\cT) \leq C_6 \bm (\cT) < \infty$. By (\ref{badmass}), $\bm ({\rm Bad} )= 0$. Then, for any compact $K \subset {\rm Bad}$, Lemma \ref{dominmass} implies that $\cP_g (K) \leq \cP^*_g (K) \leq C_6 \bm ( K) = 0$. Thus, Property {\bf Pack(2)} implies that $\cP_g ({\rm Bad})= 0$. For any Borel subset $B \subset \cT$ we get $\cP_g (B)=  \cP_g (B \cap {\rm Good} )$ and $ \bm (B)=  \bm (B \cap {\rm Good} ) $. The comparison lemma \ref{genedens} then entails  
$$ (C_\psi )^{-1} C^{-2} \bm (B\cap {\rm Good} ) \leq \cP_g (B \cap {\rm Good} )  \leq  (C_\psi )^{-1}\bm (B\cap {\rm Good} ), $$
which proves the lemma with $C_7= (C_\psi )^{-1} C^{-2}$ and $C_{8}= (C_\psi )^{-1}$. \cqfd 

\smallskip

We now prove a $0$-$1$ law to complete the proof of Theorem \ref{packingLevytreeth}. To that end we 
need the following lemma. 
\begin{lemma}
\label{technimesu} Let $\eta_n \in (0, \infty)$, $n \geq 0$, be a sequence that decreases to $0$. Then, there exists a constant $C_9 \in [0, \infty]$ that only depends on $\psi$ and on the sequence $(\eta_n)_{n \geq 0}$, such that $N$-a.e.$\; $for Lebesgue-almost all $t \in [0, \zeta]$, 
$$ \liminf_{n \rightarrow \infty} \frac{1}{\eta_n} \cP_g \big( \,  p([t,t+\eta_n]) \, \big) = C_9 \; .$$
\end{lemma}
\noi
{\bf Proof:} see Appendix \ref{appen}. \cqfd 
\begin{remark}
\label{explaintechnic} If $\cP_g (p([t, t+ \eta_n]))$ was a Borel-measurable function of $H$, 
then Lemma \ref{technimesu} would be a direct consequence of Lemma \ref{zerooneprelim}. However, this point seems difficult to prove because the packing measure $\cP_g$ is defined in two steps and the second step (\ref{packdef}) (or its variant {\bf Pack(3)}) causes measurability problems. This explains the unexpected length of the proof of Lemma \ref{technimesu}. \cq  
\end{remark}

We now completes the proof of Theorem \ref{packingLevytreeth}. 
By (\ref{massskel}), Lemma \ref{masspackequiv} and Lemma \ref{technimesu}, there exists a Borel subset $A$ of $\cC^0$ such that $N(  \cC^0 \backslash A)= 0$ and such that for any $H \in A$, the following holds true. 
\begin{description}
\item[{\bf (a)}] The mass measure $\bm$ is diffuse and $\bm ({\bf Sk} (\cT))= 0$. 
\item[{\bf (b)}] For any Borel subset $B$ of $\cT$, we have $C_7 \bm (B) \leq \cP_g (B) \leq C_8 \bm (B)$. 
\item[{\bf (c)}] For Lebesgue-almost all $t \in [0, \zeta]$, $C_9= \liminf_{n \rightarrow \infty} \eta_n^{-1} \cP_g \big( \,  p([t,t+\eta_n]) \, \big) $. 
\end{description}
We follow the end of the proof of Theorem 1.1\cite{DuLG3}. We argue deterministically on $A$. Let $0 \leq s \leq t \leq \zeta$. Suppose that there exists $r \in [0, \zeta ] \backslash [s, t]$ such that $p(r)\in p([s, t])$. Then there exists $u \in [s, t]$ such that $p(u)= p(r)$. We then get $H_u=H_r= b(u, r)$ and since $u \neq r$, (\ref{Hleafcarac}) implies that $p(r) \in {\bf Sk} (\cT)$. This easily implies the following. 
\begin{equation}
\label{intersecinterval}
p([s,t]) \cap p([u,v]) \subset  {\bf Sk} (\cT) \cup \{   p(t) \} \, , \quad 0 \leq s \leq  t \leq u \leq v \leq  \zeta \; .
\end{equation}
For any $t \in [0, \zeta]$, we set $q(t)= \cP_g (p([0, t]) )$. Observe that for any $0 \leq s \leq t \leq \zeta $, 
$p([0, t])= p([0, s]) \cup p([s, t])$. Then, (\ref{intersecinterval}) combined with ${\bf (a)}$ and ${\bf (b)}$, implies that for any $0 \leq s \leq t \leq \zeta$. 
$$ \cP_g (p([s, t]))\!=\! q(t)-q(s)  , \;  \bm (p([s, t]))\! =\! t-s,  \; {\rm and} \; C_7 (t-s) \! \leq \!q(t)-q(s) \!\leq \! C_8 (t-s). $$
The function $q$ is then absolutely continuous. Thus, $q$ is differentiable Lebesgue-almost everywhere, $q(t)-q(s)= \int_{[s, t]} q^\prime (u) du$ and for Lebesgue-almost all $t \in [0, \zeta ]$, 
$$ q^\prime (t)= \liminf_{n \rightarrow \infty} \frac{1}{\eta_n} (q(t+\eta_n)-q(t))= \liminf_{n \rightarrow \infty} \eta_n^{-1} \cP_g \big( \,  p([t,t+\eta_n]) \, \big) = C_9 \; ,  $$
by ${\bf (c)}$. This proves that $C_7 \leq C_9 \leq C_8$ and that $\cP_g (p([s, t]))= C_9 (t-s)$, for any $0 \leq s \leq t \leq \zeta$. Since $\bm$ and $\cP_g$ are diffuse, this entails that for any interval $J \subset [0, \zeta]$, $\cP_g (p(J))= C_9 \bm (p(J))= C_9 \ell (J)$, where $\ell$ stands for the Lebesgue measure.

  Let $O$ be an open set of $(\cT, d)$. Denote by $J_n$, $n \geq 0$, the (possibly empty) open connected components of $p^{-1} (O)$ that are pairwise disjoint subintervals of $[0, \zeta]$. By (\ref{intersecinterval}), ${\bf (a)}$ and ${\bf (b)}$, we get $ \cP_g \big( p(J_n) \cap p(J_{n^\prime}) \big)= 0$, for any $0 \leq n < n^\prime$. Consequently, 
\begin{eqnarray*}
 \cP_g( O) & = & \sum_{^{n \geq 0}}  \cP_g \big( p(J_n) \big)=  C_{9} \sum_{^{n \geq 0}}  
 \bm  \big( p(J_n)\big) \\
 &= & C_9 \sum_{^{n \geq 0}}  \ell  \big( J_n\big) =  C_9\ell (p^{-1} (O))= C_9 \bm (O). 
 \end{eqnarray*}
Set $c_\psi = (C_{9})^{-1}$. Then, $ c_\psi\cP_g( O)= \bm (O)$, for any open subset of $\cT$, which entails $ c_\psi\cP_g = \bm $. This completes the proof of Theorem \ref{packingLevytreeth} since $N(\cC^0 \backslash A)= 0$. \cqfd

\begin{appendix}
\section{Appendix: proof of Lemma \ref{technimesu}.}
\label{appen}
Let us recall basic facts on analytic sets and let us set some notation. Let $E$ be a topological space that is {\it Polish}. 
We shall always denote by $\cB (E)$ the Borel sigma-field of $E$. A subset $A \subset E$ is said to be {\it analytic} iff there exists an auxiliary Polish space $E^\prime$ and a Borel set $B$ of $E\times E^\prime$ equipped with the product topology such that $A$ is the projection of $B$ on $E$. We shall  denote by $\cA (E)$ the set of the analytical subsets of $E$. Then $\cB (E) \subset \cA (E)$ and 
$\cA (E)$ is stable under countable unions and intersections. Moreover, the continuous image of an analytic set is also analytic (see Jech \cite{Jech} pp 142-148 and Dudley \cite{Dudley} Chapter 13 Section 2 pp 493-501). Let $F \subset E$ be a closed subset. Since $F$ is Polish and since the canonical injection from $F$ into $E$ is continuous, we have $\cA (F) = \{ A \cap F ;  A \in \cA (E)\}$ and $\cB (F)= \{ B \cap F ; B \in \cB (E) \}$. 
We shall denote by  $\sigA (E)$ the sigma-field generated by $\cA (E)$.  
Let $\mu$ be a sigma-finite positive measure defined on $\cB(E)$. We denote by $\cB_\mu (E)$ 
the Borel sigma field $\cB (E)$ augmented by the $\mu$-negligible sets. 
Recall that $\cA (E) \subset \cB_\mu (E)$ (see Dudley \cite{Dudley} Theorem 13.2.6 p 497). This easily entails the following lemma. 
\begin{lemma}  
\label{unianaly} Let $E$ be a Polish space and let $f : E \rightarrow [0, \infty]$ be a $\sigA (E)$-measurable function. For any sigma-finite nonnegative measure $\mu$ on $\cB (E)$, there exists a $\cB(E)$-measurable function $\bar{f}_\mu : E \rightarrow [0, \infty]$ such that $\{ x \in E :\bar{f}_\mu (x)  \neq f (x) \}$ is $\mu$-negligible. 
\end{lemma}

  Recall Notation \ref{excuheightdef}: $\cC^0$ is the usual Polish space of the continuous functions from $[0, \infty)$ to $\bR$. For any $h \in \cC^0$, $\zeta(h)$ denotes $\sup \{ t \in [0, \infty): h(t) \neq 0 \}$, with the convention $\sup \emptyset = 0$. Recall that $\ccC$ stands for the set of continuous functions with compact support: $\ccC= \zeta^{-1}([0, \infty))$. For any $r \in [0, \infty)$, we also set 
$\ccC_r = \zeta^{-1} ([0, r])$ that is a closed subset of $\cC^0$. 
For any $h \in \ccC$, we set $\lVert  h \rVert = \sup |h|$. Observe that the topology of the Polish space $(\ccC_r , \lVert \cdot \rVert)$ is the trace on $\ccC_r$ of the topology of $\cC^0$. Thus, $ \cB (\ccC_r ) = \{ B \cap \ccC_r ; B \in \cB (\cC^0)\}$ and $\cA (\ccC_r ) = \{ A \cap \ccC_r ; A \in \cA (\cC^0)\}$, that are 
resp.$\; $the Borel sigma-field and the set of analytic subsets of the Polish space $(\ccC_r , \lVert \cdot \rVert)$.

  Let $h \in \ccC$. Recall from Section \ref{Levytreesec} the definition of the compact rooted $\bR$-tree $(T_h, d_h, \rho_h)$. Recall that $p_h : [0, \zeta(h)] \rightarrow T_h$ stands for the canonical projection. Recall that $\bm_h$ is the measure induced by the restriction of the Lebesgue measure $\ell$ on $[0, \zeta (h)]$ via the canonical projection (see (\ref{massmeadeter})). Recall that $\bar{p}_h : [0, \infty) \rightarrow T_h$ is given by $\bar{p}_h (t)= p_h (t \wedge \zeta (h) ) $, $t \in [0, \infty)$.

Let $g$ be the packing gauge function that is derived from $\psi$ by (\ref{gaugedef}). For any $\varepsilon \in (0, \infty)$ and for any subset $A \subset T_h$, we denote by $\cP_{^{g, h}}^{_{(\varepsilon)}} (A)$ the quantity defined by (\ref{prepremeadef}), by $\cP_{g,h}^* (A)$ the $g$-pre-packing measure of $A$ that is defined by (\ref{premeadef}). We finally denote by $\cP_{g, h} (A)$ the $g$-packing measure of $A$ that is defined by (\ref{packdef}). Property (\ref{closprepack}) asserts that $\cP_{g, h}^* (A)= \cP_{g, h}^* (\bar{A})$ where $\bar{A}$ stands for the closure of $A$ in $T_h$.  Combined with Property {\bf Pack(3)}, it entails the following: 
\begin{equation}
\label{supcompact}
\cP_{g, h} (T_h) = \inf \big\{  \sup_{n \geq 0} \cP^*_{g, h} (Q_n) \, ; \, Q_n \; {\rm compact} , \,  Q_n \subset Q_{n+1},  \, \bigcup_{^{n \geq 0}} Q_n = T_h \,  \big\} . 
\end{equation}

  Recall $C_7, C_8 \in (0, \infty) $, that appear in Lemma \ref{masspackequiv}. We introduce the following subset of functions
\begin{equation}
\label{bfSdef}
 {\bf S}= \big\{ h \in \ccC \, : \; \forall B \in \cB (T_h) \, , \;  C_7 \bm_h (B) \leq \cP_{g, h} (B) \leq C_8 \bm_h (B) \,  \big\} \; .
\end{equation} 
Lemma \ref{masspackequiv} shows that $\cC^0 \backslash {\bf S}$ is $N$-negligible. 

Recall that for any $t \in [0, \infty)$ and any $h \in \cC^0$, $\check{h}^t$ stands for the shifted function 
$(h(t+s), s \geq 0)$. Recall that for any $\eta \in [0, \infty)$ we have defined $R_\eta : \cC^0 \rightarrow \ccC$ by setting $R_\eta  h(s)= h(s\wedge \eta)-h(\eta)$, $s \in [0, \infty)$. 
Let us fix $h \in \ccC$, $t, \eta \in [0, \infty)$. To simplify notation, we set $\tilde{h}= R_\eta (\check{h}^t)$. 
Observe that for 
any $s, s^\prime \in [0, \eta]$, one has $d_h (t + s , t+ s^\prime)= d_{\tilde{h}} (s, s^\prime)$. This induces a bijective isometry $\jmath$ from $(\bar{p}_h ([t, t+ \eta] ) , d_h)$ onto $(T_{\tilde{h}} , d_{\tilde{h}})$. Moreover, $\bm_{\tilde{h}}$ is the measure induced by $\bm_h (\cdot \, \cap \bar{p}_h ([t, t+ \eta]) \,)$ via $\jmath$ and $\cP_{g, \tilde{h}}$ is the measure induced by $\cP_{g, h} (\cdot \, \cap \bar{p}_h ([t, t+ \eta]) \,)$ via $\jmath$. This entails the following lemma. 
\begin{lemma}
\label{Spiege} Let $h \in \ccC$ and let $t, \eta \in [0, \infty)$. Set $\tilde{h}= R_\eta (\check{h}^t)$. Then, 
$\cP_{g, \tilde{h}} (T_{\tilde{h}}) = \cP_{g, h} (\bar{p} ([t, t+ \eta]) )$. Moreover, if $h \in {\bf S}$, then $\tilde{h} \in {\bf S}$. 
\end{lemma}
For any $p \in \bN^*$, we denote by $\ccK_p$ the set of compact subsets of the interval $[0, p] \subset \bR$, equipped with the usual metric. We equip $\ccK_p$ with the Hausdorff distance denoted by $\dHp$. Then  $(\ccK_p , \dHp)$ is a compact metric space (see Definition \ref{Hausdorffdist} $({\bf a})$). 
\begin{lemma}
\label{mesugamma} Let $p \in \bN^*$. The function 
$ (h,K) \in  \ccC_p \times \ccK_p \mapsto  \cP_{g, h}^* (\bar{p}_h (K))\in  [0, \infty]$ is $\cB(\ccC_p) \otimes \cB (\ccK_p)$-measurable.
\end{lemma}
\noi
{\bf Proof:} For any $\varepsilon \in (0, \infty)$, $K \in \ccK_p$ and $h \in \ccC_p$, we denote by $\Pi_\varepsilon^p (h,K)$ the set of the non-empty finite sequences $(t_1 ,  r_1), \ldots , (t_n, r_n) \in  K \!\times \!(0, \varepsilon] $ such that  $d_h (t_i, t_j) >r_i+r_j$, for any $1 \leq i<j \leq n$. Since 
$\bar{p}_h (t_i)$ and $\bar{p}_h (t_j)$ can be joined by a geodesic in $T_h$, $d_h (t_i, t_j) >r_i+r_j$ is equivalent to $\bar{B}_h (\bar{p}_h (t_i) , r_i) \cap \bar{B}_h (\bar{p}_h (t_j) , r_j) = \emptyset $, for any $1 \leq i< j \leq n$. It easily entails 
$$   \cP_{g, h}^{(\varepsilon)} (\bar{p}_h (K))=  \sup \Big\{   \sum_{^{1\leq i \leq n}}\, g(r_i) \;  \, ;   \; \exists  \, (t_i ,  r_i)_{  1 \leq i \leq n } \in \Pi_\varepsilon^p (h,K) \Big\}. $$
Since $  \cP_{g, h}^* (\bar{p}_h (K))= \lim_{\varepsilon \downarrow 0} \! \downarrow \!\! \cP_{g, h}^{(\varepsilon)} (\bar{p}_h (K))$. It is sufficient to prove that for any $x \in [0, \infty)$,  the following subset $\{ (h,K) \in \ccC_p \times \ccK_p: \cP_{^{g, h}}^{_{(\varepsilon)}} (\bar{p}_h (K)) >x \}$ is open in $\ccC_p \times \ccK_p$, equipped with the product topology. So we fix  $\varepsilon, x \in (0, \infty)$, $K \in \ccK_p$ and $h \in \ccC_p$ and we assume that $\cP_{^{g, h}}^{_{(\varepsilon)}} (\bar{p}_h (K)) >x$. There exists $0 \leq t_1 < \ldots < t_n \leq p$ in $K$ and $r_1, \ldots , r_n \in (0, \varepsilon]$, such that $\sum_{1\leq i\leq n} g(r_i) >x$ and $d_h (t_i, t_j) >r_i+r_j $, for any $1 \leq i< j \leq n$. We first set 
$$ \varepsilon_1 = \min_{^{1\leq i<j\leq n}}  \!\! d_h (t_i,t_j) \!-\!r_i \!- \! r_j  \; >0 \quad {\rm and} \quad \varepsilon_2 = \frac{_1}{^3} \min_{^{1\leq i<j\leq n}} |t_i-t_j| \; >0  \; .$$
Next, set $\omega (h, \delta)= \sup \{ |h(t)-h(s)| ; s,t \in [0, \infty) : |s-t| \leq \delta \}$, for any $\delta \in (0, \infty)$ and note that $\lim_{\delta \rightarrow 0} \omega (h, \delta)= 0$. Assume that $\delta \in (0, \varepsilon_2)$ and choose $t^*_1, \ldots, t^*_n \in [0, p]$, such that $|t_i-t^*_i |< \delta $, for any $1\leq i \leq n$. An easy computation entails that 
\begin{equation} 
\label{timemove}
\big| d_h (t_i, t_j) -d_h (t^*_i, t^*_j) \big| \leq 4 \omega (h, \delta) \; , \quad 1 \leq i< j \leq n.    
\end{equation}
Next observe that for any $h^\prime \in \ccC_p$, $\big| d_h (t^*_i, t^*_j) -d_{h^\prime} (t^*_i, t^*_j) \big| \leq 4 \lVert h -h^\prime \rVert $, for any $1 \leq i< j \leq n$. Therefore, for any $\delta \in (0, \varepsilon_2)$,  for any $t^*_1, \ldots, t^*_n \in [0, p]$ such that $|t_i-t^*_i |< \delta $ and for any $h^\prime \in \ccC_p$, we get 
\begin{equation} 
\label{realmove}
\big| d_h (t_i, t_j) -d_{h^\prime} (t^*_i, t^*_j) \big| \leq 4 \omega (h, \delta) +  4 \lVert h -h^\prime \rVert \; , \quad 1 \leq i< j \leq n.    
\end{equation}
Now, fix $\delta \in (0, \varepsilon_2)$ and $\eta >0$ such that $\varepsilon_1  >4 \omega (h, \delta) + 4\eta $, which is always possible.  Let $(h^\prime , K^\prime)$ be any element of 
$\ccC_p \times \ccK_p$ such that $\lVert h-h^\prime \rVert < \eta$ and $\dHp (K, K^\prime) < \delta$.
Then, there exist $t^*_1, \ldots, t^*_n \in K^\prime$ such that $|t_i-t^*_i |< \delta $ and (\ref{realmove}) entails that $d_{h^\prime} (t^*_i, t^*_j)>r_i+r_j$, for any $1\leq i<j\leq n$. Consequently, $\cP_{g, h^\prime}^{(\varepsilon)} (\bar{p}_{h' }(K^\prime ))>x$, which completes the proof of the lemma. \cqfd 

\medskip

We next introduce the space $\ccK_p^\bN$ of the $\ccK_p$-valued sequences. We equip $\ccK_p^\bN$ with the product topology. Standard results assert that $\ccK_p^\bN$ is a Polish space (it is a compact metric space). 
We also denote by $\ccS_p$ the subset of the increasing sequences of compact subsets of $[0, p]$: 
$$ \ccS_p = \big\{ {\bf K}= (K_n)_{n \geq 0} \in \ccK_p^\bN \; : \;  K_n \subset K_{n +1} \, , \; n \geq 0 \big\} \;  \; .$$
It is easy to prove that $\ccS_p$ is a closed subset of $\ccK_p^\bN$. Therefore, $\ccS_p$ is also Polish and $\cB (\ccS_p)= \{ B \cap \ccS_p ; B \in \cB (\ccK_p)^{\otimes \bN} \}$. Recall that $\ell$ stands for the Lebesgue measure on the real line. We shall need the following Lemma. 
\begin{lemma}
\label{propricompact} The set $ \ccZ_p = \big\{ {\bf K}= (K_n )_{n \geq 0} \in \ccS_p  :  \ell \big( [0, p] \backslash \bigcup_{^{n \geq 0}} K_n \big) = 0 \big\} $ is a Borel subset of $\ccS_p$. 
\end{lemma}
\noi
{\bf Proof:}  We first prove that for any $x \in [0, \infty)$, the set $\{ K \in \ccK_p : x \leq \ell (K) \}$ is a closed subset of $(\ccK_p , \dHp)$. Observe that $\{ K \in \ccK_p : x \leq \ell (K) \} = \emptyset $ if $x >p$. We assume that $x \in [0, p]$. Let $K_n \in \ccK_p$, $n \geq 0$, be such that $x \leq \ell (K_n)$ and $\lim_n \dHp (K_n, K)= 0$. For any $\varepsilon \in (0, \infty)$, there exists $n_\varepsilon \in \bN$ such that for any $n \geq n_\varepsilon$, $K_n \subset K^{(\varepsilon)}$. Then, for any $\varepsilon \in (0, \infty)$, we have $x \leq \ell (K^{(\varepsilon)})$, which entails $x \leq \ell (K)$ by letting $\varepsilon$ go to $0$. This proves that $K \in \ccK_p \mapsto \ell (K)$ is $\cB(\ccK_p)$-measurable 
For any $ {\bf K}= (K_n )_{n \geq 0} \in \ccS_p$, we set $\Psi({\bf K})= \sup_{n \geq 0} \ell (K_n)$. 
Then, $\Psi$ is $\cB(\ccS_p)$-measurable. Now, note that $\Psi ({\bf K })= \ell (\bigcup_{n \in \geq 0} K_n)$. Thus, $\ccZ_p = \Psi^{-1} (\{ p\}) \in \cB (\ccS_p)$. \cqfd 
\begin{lemma}
\label{analyfunction} There exists a function $\Lambda : \cC^0 \rightarrow [0, \infty]$ that is $\sigA (\cC^0)$-measurable such that $\Lambda (h)= \cP_{g,h} (T_h)$, for any $h \in {\bf S}$.  
\end{lemma}
\noi
{\bf Proof:} Let us fix $p \in \bN^*$. For any ${\bf K}= (K_n)_{n \geq 0}\in \ccS_p$, and for any $h \in \ccC_p$, we set $\Gamma_p (h, {\bf K})= \sup_{n \geq 0} \cP^*_{g,h} (\bar{p}_h (K_n))$. Lemma 
\ref{mesugamma} easily implies that $\Gamma_p : \ccC_p \times \ccS_p \rightarrow [0, \infty]$ is $\cB( \ccC_p) \otimes \cB(\ccS_p) $-measurable. Recall from Lemma \ref{propricompact} the definition of $\ccZ_p$. Then, we set $ \Lambda_p (h)= \inf_{{\bf K} \in \ccZ_p} \Gamma_p ( h, {\bf K})$, for any $h \in \ccC_p$. For any $x \in (0, \infty)$, we also set 
$$ B_x =  \{ (h,{\bf K}) \in \ccC_p \times \ccS_p :   { \bf K } \in \ccZ_p \; {\rm and} \; \Gamma_p (h, {\bf K}) < x \}= (\ccC_p \times \ccZ_p ) \cap \Gamma_p^{-1} ([0, x)\, ) . $$
Note that $B_x $ is a Borel subset of $\ccC_p \times \ccS_p$. Moreover, if we denote by $\pi$ the canonical projection from $\ccC_p \times \ccS_p$ to $\ccC_p$, then $ \{ h \in \ccC_p \; : \; \Lambda_p (h) < x \} = \pi (B_x) $ that is an analytic subset of $\ccC_p$.  {\it This proves that  $\Lambda_p: \ccC_p \rightarrow [0, \infty]$ is $\sigA (\ccC_p)$-measurable}.

     We next introduce $\ccZ_p^o\!= \!\{ (K_n )_{n \geq 0} \! \in \!  \ccS_p  :   \bigcup_{^{n \geq 0}} K_n \! =\! [0, p] \}$. Observe that $\ccZ^o_p\subset \ccZ_p$.  Let us fix $h \in \ccC_p$. Note that the set of sequences of compact subsets of $T_h$ that are of the form $(\bar{p}_h (K_n) )_{n \geq 0}$, when $(K_n)_{n \geq 0}$ varies in $\ccZ_p^o$, is the same as the set of sequences of compact subsets 
$Q_n \in T_h$, $n \geq 0$, such that $Q_n \subset Q_{n+1}$ and $\bigcup Q_n =T_h$. 

Therefore, (\ref{supcompact}) entails  
\begin{equation} 
\label{geneinegpack}
\Lambda_p (h) \leq \inf_{{\bf K} \in \ccZ^o_p} \Gamma_p ( h, {\bf K} )= \cP_{g,h} (T_h) \; , \quad h \in \ccC_p .
\end{equation}
Suppose that $h \in {\bf S}\cap \ccC_p$, fix ${\bf K} = (K_n)_{n \geq 0} \in \ccZ_p$, and recall that $\cP_{g, h} (\bar{p}_h (K_n))\leq \cP^*_{g, h} (\bar{p}_h (K_n))$, for any $n \geq 0$. Then 
$$ \cP_{g,h} \big( \bar{p}_h \big( \bigcup_{^{n\geq 0}}K_n \big) \, \big) \leq \Gamma_p (h, {\bf K}) \; .$$
Since, $[0, \zeta (h)] \cap \bar{p}_h^{_{\; -1}} ( T_h \backslash  \bar{p}_h ( \bigcup K_n) \,  )  \, ) \subset [0, p] \backslash \bigcup K_n$, 
we get $ \bm_{h} \big( T_h \backslash  \bar{p}_h \big( \bigcup_{^{n\geq 0}}K_n \big) \big) =0$. Since 
$ h \in {\bf S}$, it implies  $\cP_{g, h}\big( T_h \backslash  \bar{p}_h \big( \bigcup_{^{n\geq 0}}K_n \big) \big)=0$. Consequently, $ \cP_{g, h} (T_h) \leq \Gamma_p (h, {\bf K}) $, for any ${\bf K} \in \ccZ_p$. This, combined with  (\ref{geneinegpack}), entails 
$$ \cP_{g, h} (T_h)= \Lambda_p (h)\; , \quad h \in {\bf S} \cap \ccC_p \; .$$ 
Next, it is easy to check that for any $h \in \ccC$ and for any $p, q \geq \zeta (h)$, we have $\Lambda_p (h)= \Lambda_q (h)$. Then it makes sense to set $\Lambda (h)= \Lambda_p (h)$ if $h \in \ccC$ and $p \geq \zeta (h)$ and to set $\Lambda (h)= \infty$ if $h \in \cC^0 \backslash \ccC$. Thus, for any $x \in (0, \infty)$, we get $\Lambda^{-1} ([0, x))= \bigcup \Lambda_p^{-1} ([0, x))$, 
that is an analytic subset of $\cC^0$ since $\cA (\ccC_p) \subset \cA (\cC^0)$ and since $\cA (\cC^0)$ is stable under countable unions. This completes the proof of the lemma. \cqfd 
\begin{remark}
\label{tentation} If $\ccZ^o_p $ is a Borel subset of $\ccS_p$, then the previous proof simplifies. However, we are only able to show that $ \ccS_p \backslash \ccZ^o_p $ is analytic (namely, that $\ccZ^o_p$ is co-analytic), which is not useful for our purpose. \cq 
\end{remark}
If we combine Lemma \ref{Spiege} and Lemma \ref{analyfunction}, then we get 
\begin{equation}
\label{synthese}
 \forall h \in { \bf S } \, , \; \forall t, \eta \in [0, \infty) \; , \quad \Lambda ( R_\eta \check{h}^t  )= \cP_{g, h} \big( \bar{p}_h ([t, t+ \eta] ) \, \big) \; . 
\end{equation} 
We now consider the excursion $H= (H_t)_{t \geq 0}$ of the height process. Recall notation $\zeta = \zeta (H)$. We define a measure $Q$ on $\cC^0$ as follows.  
$$ \forall B \in \cB (\cC^0)\, , \quad Q(B)= \int_0^\infty \!\!\! \!\! N( \zeta >t \, ; \,  \check{H}^t \in B) \, dt = N \Big( \int_0^\zeta \!\!\! \un_{\{ \check{H}^{t} \in B \}} dt \Big) \; . $$
For any $n \in \bN$, we set $V_n= \{ h \in \cC^0: h(0) \leq n\}$. Obviously, $V_n$ is a closed subset of $\cC^0$ and $\bigcup V_n = \cC^0$. Moreover, (\ref{meanloc}) entails $Q(V_n)\!= \! N ( \int_0^\zeta \!  \un_{\{ H_t \leq n \}} dt )\!=\! \int_0^n \!  N( L^a_\zeta) \, da \!\leq \! n $. 
{\it This proves that $Q$ is a sigma-finite measure on $\cC^0$}. 
\begin{lemma}
\label{coincanalyun}
For any $\eta \in (0, \infty)$, there exists a function $\bar{\Lambda}_\eta : \cC^0 \rightarrow [0, \infty]$ that is $\cB (\cC^0)$-measurable and such that $\{ h \in \cC^0 : \bar{\Lambda}_\eta ( R_\eta h) \neq \Lambda (R_\eta h) \}$ is $Q$-negligible. 
\end{lemma}
\noi
{\bf Proof:} We fix $\eta \in (0, \infty)$. For any $n \in \bN$, we define the finite measure $Q_{n, \eta}$ on $\cC^0$ by setting $Q_{n, \eta} (B)= Q (V_n \cap R_{\eta }^{-1} (B) \,)$, $B \in \cB(\cC^0)$. Lemma \ref{unianaly} asserts that there exists a $\cB(\cC^0)$-measurable function $\bar{\Lambda}_{n , \eta} : \cC^0 \rightarrow [0, \infty]$ such that the set $S_{n , \eta} := \{ h \in \cC^0 : \bar{\Lambda}_{n , \eta} ( h) \neq \Lambda ( h) \}$
is $Q_{n , \eta}$-negligible. Namely, $V_n \cap R_\eta^{-1} (S_{n , \eta})$ is $Q$-negligible. Let us set $S_\eta := \bigcup_{n \geq 0}V_n \cap R_\eta^{-1} (S_{n , \eta}) $ that is $Q$-negligible and let us set $\bar{ \Lambda}_\eta = \liminf_{n \rightarrow \infty} \bar{\Lambda}_{n , \eta}$ that is $\cB (\cC^0)$-measurable. It is easy to check that for any $h \in \cC^0 \backslash S_\eta$, $\bar{\Lambda}_\eta
 (R_\eta h)= \Lambda (R_\eta h)$, which completes the proof of the lemma. \cqfd 

\smallskip

We now fix a sequence $(\eta_n)_{n \geq 0}$ of positive numbers that decreases to $0$ and for any $h \in \cC^0$, we set 
$$D(h)= \liminf_{n \rightarrow \infty} \eta_{n}^{-1} \Lambda ( R_{\eta_n} h ) \quad {\rm and} \quad  \bar{D} (h)= \liminf_{n \rightarrow \infty}  \eta_{n}^{-1}\bar{\Lambda}_{\eta_n} ( R_{\eta_n} h ) \; .$$ 
Observe that $\bar{D}: \cC^0 \rightarrow [0, \infty]$ is $\cB( \cC^0)$-measurable. Lemma \ref{coincanalyun} implies that the subset $\{ h \in \cC^0: D(h) \neq \bar{D} (h) \}$ is $Q$-negligible. Moreover, Lemma \ref{zerooneprelim} entails that there exists $C_9\in [0, \infty]$ that only depends on $\psi$ and on $(\eta_n)_{n \geq 0}$ such that  
\begin{equation}
\label{zerooneborel}
Q( \{ h \in \cC^0: \bar{D} (h)\neq C_9 \})= N \Big( \int_0^\zeta \!\! \un_{\{ \bar{D} (\check{H}^t ) \neq C_9 \}}  \Big)= 0 \; .
\end{equation}
We now complete the proof of Lemma \ref{technimesu} as follows. We first set  $ {\bf S}_*= \{ h\in \cC^0 \; : \; D(h) = C_9 \} $.  Lemma \ref{coincanalyun} and (\ref{zerooneborel}) entail that $\cC_0 \backslash {\bf S}_*$ is $Q$-negligible. 
Therefore, we can find a Borel set ${\bf B}_* \in \cB(\cC^0)$ such that 
\begin{equation}
\label{Boborel}
{\bf B}_* \subset {\bf S}_* \quad {\rm and} \quad Q (\cC^0 \backslash {\bf B}_*) = 0 \; .
\end{equation}
We next set 
$$ M(h) = \int_0^{\zeta (h)} \!\!\!\! \un_{\{   \check{h}^t  \in \cC^0 \backslash {\bf S}_* \}} dt \quad {\rm and} \quad 
\bar{M}(h) = \int_0^{\zeta (h)}  \!\!\!\! \un_{\{ \check{h}^t \in \cC^0 \backslash {\bf B}_* \}} dt \; .$$
Since $\cC^0 \backslash {\bf B}_* \in \cB(\cC^0)$, standard arguments imply that $h \mapsto \bar{M} (h)$ is $\cB( \cC^0) $-measurable. Thus, the set ${\bf B}:= \{ h \in \cC^0 : \bar{M} (h)= 0 \}$ is a Borel subset of $\cC^0$. By Fubini, we get 
$$\int_{\cC^0} \bar{M} (h) N(dh) =  N \Big( \int_0^\zeta  \!\!\! \un_{\{  \check{H}^t \in  \cC^0 \backslash {\bf B}_*\}} dt \Big)= Q( \cC^0 \backslash {\bf B}_*)= 0 .$$
Therefore, $N ( \cC^0 \backslash {\bf B}) = 0 $. Recall from (\ref{bfSdef}) the definition of ${\bf S}$ and recall that $\cC^0 \backslash {\bf S}$ is $N$-negligible. Let us fix $h \in {\bf S} \cap {\bf B}$. Then, $M(h)= 0$ since $M(h) \leq \bar{M} (h)$, by (\ref{Boborel}). Namely, for Lebesgue almost all $t \in [0, \zeta(h)]$,   $ \check{h}^t  \in {\bf S}_*$, that is  $ \liminf_{n \rightarrow \infty} (\eta_n)^{-1}\Lambda ( R_{\eta_n} \check{h}^t ) = D(h) = C_9 $. 
Since $h \in {\bf S} $, (\ref{synthese})  implies that $\Lambda ( R_{\eta_n}  \check{h}^t  ) = \cP_{g, h} (\bar{p}_h ([t, t+ \eta_n ] ) )$. We thus have proved that {\it for any $h \in  {\bf S} \cap {\bf B}$, for Lebesgue-a.a.$\, t \in [0, \zeta (h) ]$, $C_9= \liminf_{n \rightarrow \infty} (\eta_n)^{-1} \cP_{g, h} (p_h ([t, t+ \eta_n ] ) )$}, which completes the proof of Lemma \ref{technimesu} since $\cC^0 \backslash ({ \bf B} \cap {\bf S})$ is $N$-negligible. \cqfd

\end{appendix}

\end{document}